\documentclass{amsart}
\usepackage{bm}
\usepackage{amssymb}
\usepackage{amsthm}
\usepackage{float}
\usepackage[shortlabels]{enumitem}
\usepackage{tikz}
\usetikzlibrary{backgrounds}
\usetikzlibrary{arrows}
\usetikzlibrary{shapes,shapes.geometric,shapes.misc}

\tikzstyle{empty dot}=[fill=none, draw=black, shape=circle]

\pgfdeclarelayer{edgelayer}
\pgfdeclarelayer{nodelayer}
\pgfsetlayers{background,edgelayer,nodelayer,main}

\tikzstyle{none}=[inner sep=0mm]

\usepackage{hyperref}
\usepackage{verbatim}
\usepackage{hyperref}
\hypersetup{
	colorlinks=true,
	linkcolor=red,
	filecolor=magenta,
	urlcolor=cyan,
	pdftitle={Overleaf Example},
	pdfpagemode=FullScreen,
}

\makeatletter
\newtheorem*{rep@theorem}{\rep@title}
\newcommand{\newreptheorem}[2]{%
	\newenvironment{rep#1}[1]{%
		\def\rep@title{#2 \ref{##1}}%
		\begin{rep@theorem}}%
		{\end{rep@theorem}}}
\makeatother

\newtheorem{theorem}{Theorem}[section]
\newreptheorem{theorem}{Theorem}

\newreptheorem{corollary}{Corollary}
\newtheorem{lemma}[theorem]{Lemma}

\newtheorem{proposition}[theorem]{Proposition}
\newreptheorem{proposition}{Proposition}

\newtheorem{claimn}{Claim}

\theoremstyle{definition}
\newtheorem{definition}[theorem]{Definition}
\newtheorem{remark}[theorem]{Remark}

\newcommand{\dlow}{\mathrm{d}}
\newcommand{\bt}[1][]{\mathsf{BT}_#1}
\newcommand{\btt}{\mathsf{BT}}
\newcommand{\btb}[1][]{\mathsf{BT}^{\mathsf{b}}_#1}
\newcommand{\sfr}{\mathsf{SFR}}

\newcommand{\bSigma}[2]{\bm{\Sigma}^#1_#2}
\newcommand{\bPi}[2]{\bm{\Pi}^#1_#2}
\newcommand{\bDelta}[2]{\bm{\Delta}^#1_#2}
\newcommand{\lSigma}[2]{\Sigma^#1_#2}
\newcommand{\lPi}[2]{\Pi^#1_#2}
\newcommand{\lDelta}[2]{\Delta^#1_#2}

\newcommand{\T}{\mathrm{T}}

\newcommand{\rca}{\ensuremath{\mathsf{RCA}_0}}
\newcommand{\wkl}{\ensuremath{\mathsf{WKL}_0}}

\DeclareMathOperator{\ran}{ran}
\newcommand{\upto}{{\upharpoonright}}

\newcommand{\N}{\ensuremath{\mathbb{N}}}

\begin{document}
	
	\title{The reverse mathematics of Brooks' Theorem}
	
\author{Alberto Marcone}
\address
	{Dipartimento di Scienze Matematiche, Informatiche e Fisiche\\
		Universit\`a di Udine\\
		33100 Udine\\
		Italy}
	\email{\href{mailto:alberto.marcone@uniud.it}{alberto.marcone@uniud.it}}
\author{Gian Marco Osso}
\address
	{Dipartimento di Scienze Matematiche, Informatiche e Fisiche\\
		Universit\`a di Udine\\
		33100 Udine\\
		Italy}
	\email{\href{mailto:osso.gianmarco@spes.uniud.it}{osso.gianmarco@spes.uniud.it}}
\date{\today}
\thanks{Both authors were supported by the Italian PRIN 2022 ``Models, sets and classifications'', prot.\ 2022TECZJA.
Marcone is a member of INdAM-GNSAGA}
	
	
	\begin{abstract}
		This is an analysis of the status of Brooks' Theorem, a celebrated result in graph coloring, from the point of view of Reverse Mathematics. We prove that the restriction of Brooks' theorem to bounded graphs of degree greater than or equal to $3$ is provable in $\rca$, while the statement for arbitrary graphs is equivalent to $\wkl$ over $\rca$. Brooks' Theorem for degree $2$, even when restricted to bounded graphs, is equivalent to $\wkl$ over $\rca$.
	\end{abstract}
	
	\maketitle
		
	\section{Introduction}
	
	This paper is a contribution to the reverse mathematics program, which seeks to identify the axioms necessary to prove theorems of ordinary mathematics. We analyze Brooks' theorem, an old result about graph colorings which provides an upper bound on the chromatic numbers of graphs of finite degree (i.e.\ the graphs for which there is some $m \in \N$ such that all nodes have at most $m$ neighbors). If $m$ is the maximum degree of a graph $G$, the greedy coloring algorithm yields a coloring in $m+1$ vertices. Brooks' Theorem (\cite{Brooks}) identifies exactly the class of graphs for which this coloring turns out to be optimal (meaning that any other proper coloring uses at least $m+1$ colors).
	
	\begin{theorem}
		Let $d \geq 2$ and let $G$ be a graph such that all nodes in $G$ have at most $d$ neighbors. Suppose that $G$ does not have any $d+1$-clique as a subgraph. If $d=2$, additionally assume that $G$ does not have any odd cycles. Then $G$ is $d$-colorable.
	\end{theorem}
	
	This result is a classic of graph coloring, and a source of inspiration for combinatorialists to this day (see \cite{btab}).
	
	A very important result in (noncomputable) graph colorings due to De Bruijn and Erd\H{o}s, which will appear throughout the paper, is the following:
	
	\begin{theorem}\label{thm:comp}
		Let $G$ be an infinite graph and $k \geq 2$, then $G$ is $k$-colorable if and only if every finite subgraph of $G$ is $k$-colorable.
	\end{theorem}
	
	Theorem \ref{thm:comp} can be proven using the compactness theorem for first order logic. In the context of reverse mathematics, the compactness theorem for first order logic is equivalent to $\wkl$ (see \cite[Theorem IV.3.3]{simpson}). In fact, it is known that the De Bruijn-Erd\H{o}s theorem is also equivalent to $\wkl$ (see, e.g.\ \cite[Theorem 2]{hirstgasarch}).
	
	We will use Theorem \ref{thm:comp} to show that, unsurprisingly, Brooks' Theorem is provable in $\wkl$. Several subcases of Brooks' Theorem are shown equivalent to $\wkl$. These are the only (to our knowledge) examples of reversals of statements in graph coloring to $\wkl$ which only employ graphs of finite degree. Other subcases are actually provable in $\rca$. We show this latter fact by clarifying an old proof of Tverberg \cite{Tverberg84}, and showing that it is formalizable in $\rca$. A distinction which is relevant in our study is that between graphs of finite degrees and bounded graphs (this corresponds to the distinction between computable graphs and so-called ``highly recursive'' graphs).
	
	All our proofs are completely uniform, so that they correspond naturally to Weihrauch reductions. In particular, the versions of Brooks' theorem which we show provable in $\rca$ are, in the context of Weihrauch reductions, reducible to the identity function. Similarly, the variants equivalent to $\wkl$ correspond to problems Weihrauch equivalent to the natural correspondent of $\wkl$, denoted $\mathsf{C}_{2^{\mathbb{N}}}$ in the literature on Weihrauch degrees (see \cite[Section 11]{BGP}).
	
	\section{Preliminaries}
	
	\subsection{Basic graph theoretic notions} We introduce the terminology used to deal with graphs and graph colorings.
	
	\begin{definition}
		A \emph{graph} $G$ is a pair $(V,E)$ where $V$ is a nonempty set (the \emph{vertices} of $G$) and $E$ is a \emph{symmetric} and \emph{anti-reflexive} ($(x,x) \notin E$ for every $x$) subset of $V^2$. $E$ is the set of \emph{edges} of $G$.
	\end{definition}
	A pair of vertices $x,y \in V$ is said to be \emph{adjacent} (or connected via an edge) if $(x,y) \in E$. Given a vertex $x$, a vertex $y$ which is adjacent to $x$ is called a \emph{neighbor} of $x$. The set of neighbors of $x$ is denoted $N(x)$.
	
	\begin{definition}
		The number of neighbors of a vertex $x$ is called the \emph{degree} of of $x$, denoted as $\delta(x)$.
	\end{definition}
	A graph is called \emph{locally finite} if every vertex has finite degree.
	
	\begin{definition}
		A \emph{graph of finite degree} is a locally finite graph $G=(V,E)$ such that there is $m \in \N$ such that, for every $x \in V$, $\delta(x) \leq m$. If $G=(V,E)$ is a graph of finite degree and $\Delta(G)=\max\{\delta(x): x \in V\}$ we say that $\Delta(G)$ is the \emph{degree} of $G$. A graph where each vertex has exactly $d$ neighbors is called \emph{$d$-regular}.
	\end{definition}
	
	We give definitions of neighborhoods, induced subgraphs, cycles and cliques.
	
	\begin{definition}
		Let $G=(V,E)$ be a graph and let $x \in V$. The \emph{$1$-neighborhood} of $x$, denoted $N_1(x)$ is the set $\{x\} \cup \{y \in V: xEy\}$. The $1$-neighborhood of a set $U \subseteq V$ is the set $\bigcup_{x \in U} N_1(x)$. For any $n \geq 1$ and any $U \subseteq V$, the $n+1$-neighborhood of $U$, denoted $N_{n+1}(U)$ is the set $N_1(N_n(U))$.
	\end{definition}
	
	Notice that $N(x)$ and $N_1(x)$ denote different sets: indeed, $x$ never belongs to the former and always belongs to the latter.
	
	\begin{definition}
		Let $G=(V,E)$ be a graph. An \emph{induced subgraph} of $G$ is a graph $(U,F)$ such that $U \subseteq V$ and $F = U^2 \cap E$.
	\end{definition}
	
	Given a graph $(V,E)$ and a set of vertices $U \subseteq V$ we often silently identify $U$ with the induced subgraph $(U,U^2 \cap E)$. All mentions of ``subgraphs'' in the rest of this paper refer to induced subgraphs.
	
	\begin{definition}
		For any $n \in \N$ an \emph{$n$-clique} (or \emph{clique on $n$ vertices}) is a graph $(V,E)$ where $|V|=n$ and for all $v \neq w \in V$, $(v,w) \in E$.
	\end{definition}
	
	In other words, an $n$-clique is a graph consisting of $n$ vertices in which all pairs of vertices are connected by an edge.
	
	\begin{definition}
		Let $G=(V,E)$ be a graph and $n \in \N$. A \emph{path of length $n$} in $G$ is a sequence $(v(i) : i < n)$ of pairwise distinct vertices such that $v(i)Ev(i+1)$ for all $i<n-1$. If $P=(v(i) : i < n)$ is a path, we say that $P$ is a \emph{path from $v(0)$ to $v(n-1)$}. If $(v(i) : i < n)$ is a path and $v(0) E v(n-1)$, then the sequence $(v(i) : i \leq n)$ where $v(n)=v(0)$ is a \emph{cycle} in $G$. A cycle is called \emph{odd} if it contains an odd number of vertices, \emph{even} if it contains an even number of vertices.
	\end{definition}
	
	We also define the concept of connected component.
	
	\begin{definition}
		Let $G=(V,E)$ be a graph and let $x$ and $y$ be vertices of $G$. We say that $x$ and $y$ are \emph{connected} if there exists a path from $x$ to $y$. Connectedness in a graph is an equivalence relation on the vertex set $V$. Equivalence classes of the connectedness relations are called \emph{connected components} (or just \emph{components}); we denote by $C_x$ the connected component containing $x$. A graph is \emph{connected} if it consists of a single connected component.
	\end{definition}
	
	The following kinds of graphs play a role in the paper.
	
	\begin{definition}
		Let $G=(V,E)$ be a connected graph with $\Delta(G)=2$, then $G$ is a \emph{circle} if $V$ is finite and all nodes have degree $2$. If $V$ is finite but some nodes have degree less than $2$, then there are exactly two nodes (the \emph{endpoints} of $G$) with degree $1$. In this case, we say that $G$ is a \emph{line segment}.
	\end{definition}
	
	We now define the concepts of a proper graph coloring and chromatic number.
	
	\begin{definition}
		Let $G=(V,E)$ be a graph and $X$ be any set. A \emph{proper $X$-coloring} of $G$ is a function $c \colon V \rightarrow X$ such that for every $x,y \in V$, $xEy$ implies that $c(x) \neq c(y)$. If $X=\{0, \dots, n-1\}$, we say that $c$ is a proper $n$-coloring.
	\end{definition}
	\begin{definition}
		Let $G$ be a graph such that there exists a proper coloring of $G$ using finitely many colors. The least $n$ for which there is a proper $n$-coloring of $G$ is the \emph{chromatic number} of $G$, denoted $\chi(G)$.
	\end{definition}
	
	A greedy coloring algorithm readily shows that if $G$ is a graph of finite degree, then $\chi(G) \leq \Delta(G)+1$.
	
	\subsection{Graphs in $\rca$}
	
	We now give definitions necessary to deal with graphs in $\rca$. These are motivated by considerations of effectiveness, and, as usual, have counterparts in computability theory.
	
	\begin{definition}
		A \emph{bounded} graph is a triple $(V,E,f)$ where $(V,E)$ is a locally finite graph and $f \colon V \rightarrow \N$ maps each $x \in V$ to its degree $\delta(x)$.
	\end{definition}
	
	A key computational difference between locally finite graphs and bounded graphs is that, in a bounded graph, we can compute neighborhoods of a given point (or set of points), whereas this is not always possible in arbitrary locally finite graphs. Indeed, we have:
	
	\begin{remark}
		If $G=(V,E,f)$ is a bounded graph with $V  \subseteq \omega$, then the function $(x,n) \mapsto N_n(x)$ is $G$-computable.
		
		Conversely, if $G=(V,E)$ is a locally finite graph with $V \subseteq \omega$ and such that the function $x \mapsto N_1(x)$ is computable, then there is a $G$-computable $f$ such that $f(x)=\delta_G(x)$ (so, $G$ is computable as a bounded graph).
		
		In terms of computability, if $G=(V,E)$ is an arbitrary graph, then for every $n \geq 1$ and every set $U \subseteq V$, $N_n(U)$ is $\lSigma{0}{1}$-definable from $U \oplus G$.
		
		On the other hand, if $G=(V,E,f)$ is a bounded graph and $U \subseteq V$, then the formalization of the algorithm mentioned above shows that the sequence $(N_i(U): i \geq 1)$ is $U \oplus G$-computable, and hence $\rca$ proves its existence as a set.
	\end{remark}
	
	In the previous section, we gave a definition of connectedness and connected components. It is not hard to see that for every countable graph $G=(V,E)$ and every node $x \in V$, $C_x$ is uniformly $\lSigma{0}{1}$-definable in $G$ and $x$. Moreover, it is known that there is a computable bounded graph $G=(V,E,f)$ and a node $x \in V$ such that $C_x \equiv_{\T} \emptyset'$, so that we cannot hope for a more effective definition.
	It follows that $\rca$ does not prove the existence of connected components in bounded graphs (see \cite[Theorem 2.1 and Corollary 2.2]{GHM}). Nonetheless, we can always recognize a connected component of bounded graphs when we see one:
	
	\begin{remark}\label{rem:connected}
		Let $G=(V,E,f)$ be a bounded graph and let $H \subseteq V$ be a set of vertices. We say that $H$ is a connected component of $G$ if and only if $N_1(H)=H$. This is clearly expressible in $\rca$. Moreover, if $H$ is finite, the formula ``$H$ is a connected component of $G$'' is $\bDelta{0}{1}$.
	\end{remark}
	
	In the rest of the paper we use the expression $\chi(G) \leq d$ as an abbreviation for ``$G$ admits a proper $d$-coloring''.
	Indeed, it is not clear whether the notion of chromatic number is well defined in $\rca$: the statement ``any graph $G$ which admits a proper coloring in finitely many colors also has a chromatic number'' follows from the least number principle for $\bSigma{1}{1}$ classes, and it can be seen to also follow from $\wkl$, via an application of Theorem \ref{thm:comp}.
	
	We conclude the section with the formal statement of Brooks's Theorem. More precisely, we give one statement for each possible value of $\Delta(G)$. We keep the case $\Delta(G)=2$ separate, as it has stronger hypotheses.
	
	\begin{definition}
		We define $\bt 2$ as (the formalization in second order arithmetic of) the statement ``every graph $G$ with $\Delta(G)=2$ and without odd cycles is such that $\chi(G) \leq 2$''.
	\end{definition}
	
	\begin{definition}
		We define $\bt d$ as (the formalization in second order arithmetic of) the statement ``every graph $G$ with $\Delta(G)=d$ which induces no $d+1$ cliques is such that $\chi(G) \leq d$''. For every $d \geq 2$, we also denote by $\btb d$ the restriction of $\bt d$ to bounded graphs.
	\end{definition}
	
	We remark that by $\bt d$ and $\btb d$ we intend formulae in the language of first order arithmetic with a free number variable $d$. Thus, $d$ is not restricted to being a standard natural number.
	
	\begin{definition}
		We define $\btt$, which we regard as the correct second order arithmetic formulation of Brooks' Theorem, as $\bt 2 \land \forall d \geq 3\,\, \bt d$.
	\end{definition}
	
	\section{The strength of the full Brooks' Theorem}
	
	As one would expect, Brooks' Theorem is provable in $\wkl$. This is due to the fact that Brooks' theorem for finite graphs is provable in $\rca$, and the compactness argument behind the proof of the De Bruijn-Erd\H{o}s theorem goes through in $\wkl$. We now prove these facts.
	
	\begin{lemma}\label{lem:brooksfinitebi}
		The restriction of $\bt 2$ to bounded graphs with no infinite connected component is provable in $\rca$. In particular, the restriction of $\bt 2$ to finite graphs is provable in $\rca$.
	\end{lemma}
	\begin{proof}
		Work in $\rca$. Let $G=(V,E,f)$ be a bounded graph with $\Delta(G)=2$, no odd cycles, and no infinite connected component.		
		We build our coloring exploiting the fact that the connected components of $G$ are either singletons or line segments or circles with an even number of vertices.
		Since every connected component is finite we have a listing $(C(i) : i \in \N)$ of all the connected components of $G$, and we can color them one by one.
		
		Given $i$, if $C(i)$ is a singleton $\{x\}$, we color $x$ with color $0$. If $C(i)$ contains two nodes $y < z$ with $f(y)=f(z)=1$ (this occurs if $C(i)$ is a line segment with endpoints $y$ and $z$), then we $2$-color the vertices along the path from $y$ to $z$, starting from $y$. Otherwise, $C(i)$ is a circle and we $2$-color the vertices therein starting by giving color $0$ to the least node of $C(i)$, and coloring each successive vertex with the only acceptable color. The fact that the circle has an even number of vertices guarantees that the resulting coloring is proper.
		
		This construction is computable and guaranteed to result in a $2$-coloring of $G$. This shows that the restriction of $\bt 2$ to bounded graphs with no infinite connected component is provable in $\rca$.
		
		To show that the restriction of $\bt 2$ to finite graphs is provable in $\rca$, notice that $\rca$ proves that every finite graph is bounded.
	\end{proof}
	\begin{lemma}\label{lem:brooksfinited}
		The restriction of $\forall d \geq 3 \,\, \bt d$ to finite graphs is provable in $\rca$.
	\end{lemma}
	\begin{proof}
		We can follow Brooks' proof of his theorem for finite graphs in \cite{Brooks}, and we see that the proof is essentially an application of $\bDelta{0}{1}$-induction, a principle which is allowed in $\rca$.
	\end{proof}
	
	\begin{theorem}\label{thm:upperbound}
		Brooks' Theorem $\btt$ is provable in $\wkl$.
	\end{theorem}
	\begin{proof}
		Let $d \in \N$ and let $G=(V,E)$ be a graph with $\Delta(G)=d$. If $d=2$, assume that $G$ does not induce odd cycles, and if $d \geq 3$ assume that $G$ does not induce $d+1$ cliques.
		
		Since the restriction of $\btt$ to finite graphs is provable in $\rca$, (by Lemma \ref{lem:brooksfinitebi} for the case $d=2$ and Lemma \ref{lem:brooksfinited} for the other cases) it follows that any finite subgraph of $G$ is $d$-colorable. By the De Bruijn-Erd\H{o}s Theorem, which is provable in $\wkl$, it follows that $G$ is $d$-colorable.
	\end{proof}
	
	We now show that $\btt$ is equivalent to $\wkl$ over $\rca$, providing a reversal of Theorem \ref{thm:upperbound}. Indeed, we show that even restricted forms of Brooks' Theorem are already equivalent to $\wkl$ over $\rca$.
	
	We introduce the following well-known principle.
	\begin{definition}
		This definition is made within $\rca$. Define the principle $\sfr$ (for \emph{Separating Function Ranges}) as (the formalization in second order arithmetic) of the statement ``for every pair of injective functions $f \colon \N \rightarrow \N$ and $g \colon \N \rightarrow \N$ if $f(n) \neq g(m)$ for all $n,m \in \N$, then there is a set $X$ such that for all $n \in \N$, $f(n) \in X$ and $g(n) \notin X$''.
	\end{definition}
	
	The following Proposition is Lemma IV.4.4 in \cite{simpson}.
	
	\begin{proposition}
		The principles $\wkl$ and $\sfr$ are equivalent over $\rca$.
	\end{proposition}
	
	We now prove the main theorem of the section. The key construction of the proof (which also features in the proof of Proposition \ref{prop:boundedbiwkl}) was vaguely inspired by Bean's construction of a $3$-colorable computable graph with no computable finite coloring (see \cite{bean}).
	
	\begin{theorem}\label{thm:mainrev}
		For any $d \geq 2$, we have that $\bt d$ implies $\wkl$ over $\rca$.
	\end{theorem}
	\begin{proof}
		First we give a proof that $\bt 2$ implies $\sfr$ over $\rca$. This construction is the blueprint for the constructions showing that $\bt d$ implies $\wkl$ for any $d \geq 3$.
		
		So, suppose $f$ and $g$ are injective functions from $\N$ to $\N$ with disjoint ranges. We define a graph $H_{2}=(V_{2},E_{2})$ as follows: $V_{2}=\{a_n, b_n, c_n, c'_n, d'_n: n \in \N\}$ and
\begin{multline*}
E_2=\{(a_n,c_m) , (b_n, c_m) : f(n)=m\} \cup \{(a_n,c'_m), (b_n, d'_m): g(n)=m\} \cup {}\\
\{(c'_m,d'_m) : m \in \N\}.
\end{multline*}

Intuitively, we have a pair of vertices $a_n$ and $b_n$ for every natural
number. We will use the colors of these vertices to separate $\ran(f)$ from
$\ran(g)$. Moreover, we have auxiliary vertices $c_m$, $c'_m$ and $d'_m$ for
every $m$. Edges of $H_2$ are as detailed in Figure \ref{fig2}. If $n$ is not
in $\ran(g) \cup \ran(f)$, both $a_n$ and $b_n$ are isolated vertices of $G$.
		\begin{figure}
			\centering
			{\begin{tikzpicture}
					\begin{pgfonlayer}{nodelayer}
						\node [style=empty dot] (0) at (-2, 1) {$a_n$};
						\node [style=empty dot] (1) at (-2, -1) {$c_m$};
						\node [style=empty dot] (2) at (0, 1) {$b_n$};
						\node [style=empty dot] (4) at (1.5, 1) {$a_k$};
						\node [style=empty dot] (5) at (1.5, -1) {$c'_h$};
						\node [style=empty dot] (6) at (3.5, 1) {$b_k$};
						\node [style=empty dot] (7) at (3.5, -1) {$d'_h$};
						\node [style=empty dot] (8) at (5, 1) {$a_i$};
						\node [style=empty dot] (9) at (7, 1) {$b_i$};
					\end{pgfonlayer}
					\begin{pgfonlayer}{edgelayer}
						\draw (0) to (1);
						\draw (1) to (2);
						\draw [in=90, out=-90] (4) to (5);
						\draw (5) to (7);
						\draw (7) to (6);
					\end{pgfonlayer}
				\end{tikzpicture}
				
			}	
			\caption{\label{fig2}A portion of the graph $H_{2}$ for a pair of functions $f$ and $g$ such that $f(m) = n$, $g(h) = k$, and $i \notin \ran(f) \cup \ran(g)$.}
		\end{figure}
		Since $H_{2}$ is $\lDelta{0}{1}$-definable from $f$ and $g$, it exists in $\rca$. Moreover, it is clear from the definition that $\Delta(H_{2})=2$ and $H_{2}$ does not induce any (odd) cycles. So, by $\bt 2$, there is a $2$-coloring $c$ of $H_{2}$. It is easy to see that if $f(m)=n$, it must be that $c(a_n)=c(b_n)$, while if $g(m)=n$ it must be that $c(a_n) \neq c(b_n)$. So we can define $X=\{n \in \N: c(a_n)=c(b_n)\}$. The set $X$ is $\lDelta{0}{1}$-definable from $c$, so it exists in $\rca$. Moreover, we have that if $n \in \ran(f)$, then $n \in X$, and if $n \in \ran(g)$, then $n \notin X$. This shows that $\sfr$ holds.
		
		The proof that for any $d \in \N$, $\bt d$ implies $\sfr$ is entirely analogous. For a fixed $d$, given $f$ and $g$ functions with disjoint ranges, we can define a graph $H_{d}$ as follows. We start with a pair of vertices $a_n$ and $b_n$ for each natural number $n$. We then have, for every $m \in \omega$, two copies $G_m$ and $G'_m$ of $K_{d-1}$, and an auxiliary node $e'_m$ which is connected to every node in $G'_m$. If $f(m)=n$, we add an edge between $a_n$ and every node of $G_m$, and an edge between $b_n$ and every node of $G_m$. If $g(m)=n$, we add an edge between $a_n$ and every node of $G'_m$, and an edge between $e'_m$ and $b_n$. We have that, if $f(m)=n$, then the subgraphs $\{a_n\} \cup G_m$ and $\{b_n\} \cup G_m$ are both $d$-cliques. This readily implies that any proper $d$-coloring of $H_d$ assigns the same color to both $a_n$ and $b_n$. On the other hand, if $g(m)=n$, the same reasoning shows that any proper $d$-coloring of $H_d$ assigns the same color to $a_n$ and $e'_m$, and hence the color assigned to $b_n$ differs from that assigned to $a_n$. The construction for $d=3$ is exemplified in Figure \ref{fig3}.
		\begin{figure}
			\centering
			\begin{tikzpicture}
				\begin{pgfonlayer}{nodelayer}
					\node [style=empty dot] (0) at (-4, 0) {$a_n$};
					\node [style=empty dot] (1) at (-3, -1) {$d_m$};
					\node [style=empty dot] (2) at (-2, 0) {$b_n$};
					\node [style=empty dot] (3) at (-3, 1) {$c_m$};
					\node [style=empty dot] (10) at (-1, 0) {$a_k$};
					\node [style=empty dot] (11) at (0, -1) {$d'_h$};
					\node [style=empty dot] (12) at (1, 0) {$e'_h$};
					\node [style=empty dot] (13) at (0, 1) {$c'_h$};
					\node [style=empty dot] (14) at (2, 0) {$b_k$};
					\node [style=empty dot] (15) at (3, 0) {$a_i$};
					\node [style=empty dot] (16) at (5, 0) {$b_i$};
				\end{pgfonlayer}
				\begin{pgfonlayer}{edgelayer}
					\draw (0) to (1);
					\draw (1) to (2);
					\draw (2) to (3);
					\draw (3) to (0);
					\draw (3) to (1);
					\draw (10) to (11);
					\draw (11) to (12);
					\draw (12) to (13);
					\draw (13) to (10);
					\draw (13) to (11);
					\draw (12) to (14);
				\end{pgfonlayer}
			\end{tikzpicture}
			\caption{\label{fig3} A portion of the graph $H_{3}$ for a pair of functions $f$ and $g$ such that $f(m) = n$, $g(h) = k$, and $i \notin \ran(f) \cup \ran(g)$.}
		\end{figure}
		
		Again we have that $H_{d}$ is $\lDelta{0}{1}$-definable from $f$ and $g$ and moreover $\Delta(H_{d}) = d$ as it is clear that no vertex has more than $d$ neighbors and there exists vertices with exactly $d$ neighbors (in the case $d=3$ pictured in Figure \ref{fig3}, the vertices $e'_m$ will have exactly $3$ neighbors). It is also evident that $H_{d}$ does not induce any $d+1$-clique. Exactly as in the proof for $d=2$, we can see that if $c$ is a $d$-coloring of $H_{d}$, then $X=\{n \in \N : c(a_n)=c(b_n)\}$ separates the range of $f$ from the range of $g$. Hence, $\bt d$ implies $\sfr$.
	\end{proof}
	\begin{remark}
		The graphs used in the proof of Proposition \ref{prop:boundedbiwkl} have no infinite connected component. This shows that the boundedness assumption in Lemma \ref{lem:brooksfinitebi} is necessary. By Lemma \ref{lem:brooksfinitebi} and Theorem \ref{thm:mainrev}, we infer that, in general, $\rca$ does not prove the existence of the degree function for $H_{2}$ (in particular in $\rca$ does not show that the set of isolated node in $H_2$ exists). It is also easy to see that $\rca$ does not in general prove the existence of the degree functions for $H_{d}$ for any $d \geq 3$.
	\end{remark}
	
	We summarize the results of the section.
	
	\begin{theorem}
		Over $\rca$, for every $d \geq 2$, the following are equivalent:
		\begin{itemize}
			\item $\wkl$,
			\item $\bt d$,
			\item $\btt$.
		\end{itemize}
	\end{theorem}

	\section{Brooks' Theorem for bounded graphs}\label{sec:bounded}
	
	We now analyze the strength of Brooks' Theorem restricted to bounded graphs. It is well known that the computational properties of bounded graphs differ from those of graphs of finite degree whose degree function is not given (see for example \cite{hirstgasarch}). This difference is visible in our setting too.
	
	In particular, one tool that we will use throughout the section is the following (this is a remark in both \cite{Schmerl82} and \cite{Tverberg84}):
	
	\begin{lemma}\label{lem:boundedregular}
		The following is provable in $\rca$. Let $G=(V,E,f)$ be a bounded graph with $\Delta(G)=d$ and no $d+1$-clique. There exists a $d$-regular graph $G'=(V',E')$ which has no $d+1$-cliques, a subset $U' \subseteq V'$ and a function $h \colon V \rightarrow U'$ which induces an isomorphism between $G$ and the induced subgraph of $G'$ determined by $U'$.
	\end{lemma}
	\begin{proof}
		Let $G^0=G$. Recursively define, for $0 \leq i<d$, the graph $G^{i+1}=(V_{i+1}, E_{i+1}, f_{i+1})$ as the result of taking two disjoint copies of $G^i$, say $G^i_0$ and $G^i_1$ and adding an edge between any $x_0 \in G^i_0$ with $f_{i}(x_0)<d$ and the corresponding $x_1$ in $G^i_1$. Now let $x \in V$ be such that $\delta(x)=k \leq d$. By $\bDelta{0}{1}$-induction one can see that for all $0<i\leq d$, all copies of $x$ in $G^i$ have degree $\min\{k+i,d\}$. In particular, all copies of $x$ have degree $d$ in $G^d$, so $G^d$ is $d$-regular. It is clear by construction that $G^d$ has no $d+1$-clique and it is easy to define $U'$ and $h$.
	\end{proof}
	
	Since it is clear that a $d$-coloring of $G'$, together with the function $h$, allows us to recover a $d$-coloring for $G$, we have that, for every $d \geq 2$, $\btb[d]$ is equivalent to the restriction of $\btb[d]$ to $d$-regular graphs. We use this in the proofs of Theorem \ref{thm:bt3}, Proposition \ref{prop:ddown}, and Theorem \ref{thm:btd}.
	
	\subsection{The strength of $\btb[2]$}
	
	We show that $\btb[2]$ is equivalent to $\wkl$ over $\rca$. This in particular implies that the well known characterization of bipartite graphs as ``graphs with no odd cycles'' is already at the level of $\wkl$ even for bounded graphs of degree $2$. The fact that the characterization of arbitrary bipartite graphs is at the level of $\wkl$ is an old result of Hirst, see \cite[Theorem 3.3]{mtrm}.
	
	\begin{proposition}\label{prop:boundedbiwkl}
		The principle $\btb 2$ is equivalent to $\wkl$ over $\rca$.
	\end{proposition}
	\begin{proof}
		As we mentioned above, we show that the restriction of $\bt 2$ to bounded graphs implies $\sfr$.
		
		Let $f \colon \N \rightarrow \N$ and $g \colon \N \rightarrow \N$ be two injective functions such that for every $n$ and $m$, $f(n) \neq g(m)$. We build a bounded graph $G=(V,E)$ (depicted in Figure \ref{fig:bt2}) as follows: the vertices $V$ are of the form $\{l_{n,k}, m_{n,k}, r_{n,k}: n \in \N, k \in \N\}$.
		We add edges on the basis of the functions $f$ and $g$. For every $n$ and every $k$, if $f(m) \neq n$ and $g(m) \neq n$ for every $m \leq k$, then $l_{n,k}\,E\,l_{n,k+1}$ and $r_{n,k}\, E\, r_{n,k+1}$. If $k$ is such that $f(k)=n$, then we set $l_{n,k} \, E\, r_{n,k}$, while if $k$ is such that $g(k)=n$, we set $l_{n,k} \, E\, m_{n,k}$ and $m_{n,k} \, E\, r_{n,k}$. We add no other edge.
		
		\begin{figure}
			\centering
			\scalebox{0.82}{\begin{tikzpicture}
					\begin{pgfonlayer}{nodelayer}
						\node [style=empty dot] (0) at (-7, -5.5) {$l_{1,0}$};
						\node [style=empty dot] (1) at (-7, -4) {$l_{1,1}$};
						\node [style=empty dot] (2) at (-7, -2.5) {$l_{1,2}$};
						\node [style=empty dot] (3) at (-7, -1) {$l_{1,3}$};
						\node [style=empty dot] (4) at (-5.5, -1) {$m_{1,3}$};
						\node [style=empty dot] (5) at (-4, -1) {$r_{1,3}$};
						\node [style=empty dot] (6) at (-4, -2.5) {$r_{1,2}$};
						\node [style=empty dot] (7) at (-4, -4) {$r_{1,1}$};
						\node [style=empty dot] (8) at (-4, -5.5) {$r_{1,0}$};
						\node [style=empty dot] (9) at (-5.5, -5.5) {$m_{1,0}$};
						\node [style=empty dot] (10) at (-5.5, -4) {$m_{1,1}$};
						\node [style=empty dot] (11) at (-5.5, -2.5) {$m_{1,2}$};
						\node [style=empty dot] (12) at (-7, 0.5) {$l_{1,4}$};
						\node [style=empty dot] (15) at (-5.5, 0.5) {$m_{1,4}$};
						\node [style=empty dot] (18) at (-4, 0.5) {$r_{1,4}$};
						\node [style=empty dot] (21) at (-1, -5.5) {$l_{2,0}$};
						\node [style=empty dot] (22) at (-1, -4) {$l_{2,1}$};
						\node [style=empty dot] (23) at (-1, -2.5) {$l_{2,2}$};
						\node [style=empty dot] (24) at (-1, -1) {$l_{2,3}$};
						\node [style=empty dot] (25) at (-1, 0.5) {$l_{2,4}$};
						\node [style=empty dot] (30) at (0.5, 0.5) {$m_{2,4}$};
						\node [style=empty dot] (31) at (0.5, -1) {$m_{2,3}$};
						\node [style=empty dot] (32) at (0.5, -2.5) {$m_{2,2}$};
						\node [style=empty dot] (33) at (0.5, -4) {$m_{2,1}$};
						\node [style=empty dot] (34) at (0.5, -5.5) {$m_{2,0}$};
						\node [style=empty dot] (35) at (2, -5.5) {$r_{2,0}$};
						\node [style=empty dot] (36) at (2, -4) {$r_{2,1}$};
						\node [style=empty dot] (37) at (2, -2.5) {$r_{2,2}$};
						\node [style=empty dot] (38) at (2, -1) {$r_{2,3}$};
						\node [style=empty dot] (39) at (2, 0.5) {$r_{2,4}$};
						\node [style=empty dot] (43) at (-13, -5.5) {$l_{0,0}$};
						\node [style=empty dot] (44) at (-11.5, -5.5) {$m_{0,0}$};
						\node [style=empty dot] (45) at (-10, -5.5) {$r_{0,0}$};
						\node [style=empty dot] (46) at (-13, -4) {$l_{0,1}$};
						\node [style=empty dot] (47) at (-13, -2.5) {$l_{0,2}$};
						\node [style=empty dot] (48) at (-13, -1) {$l_{0,3}$};
						\node [style=empty dot] (49) at (-11.5, -1) {$m_{0,3}$};
						\node [style=empty dot] (50) at (-11.5, -2.5) {$m_{0,2}$};
						\node [style=empty dot] (51) at (-11.5, -4) {$m_{0,1}$};
						\node [style=empty dot] (52) at (-10, -4) {$r_{0,1}$};
						\node [style=empty dot] (53) at (-10, -2.5) {$r_{0,2}$};
						\node [style=empty dot] (54) at (-10, -1) {$r_{0,3}$};
						\node [style=empty dot] (55) at (-10, 0.5) {$r_{0,4}$};
						\node [style=empty dot] (56) at (-11.5, 0.5) {$m_{0,4}$};
						\node [style=empty dot] (57) at (-13, 0.5) {$l_{0,4}$};
						\node [style=none] (63) at (-10, 2) {};
						\node [style=none] (64) at (-13, 2) {};
					\end{pgfonlayer}
					\begin{pgfonlayer}{edgelayer}
						\draw (0) to (1);
						\draw (1) to (2);
						\draw (2) to (3);
						\draw (5) to (6);
						\draw [in=90, out=-90] (6) to (7);
						\draw (7) to (8);
						\draw (21) to (22);
						\draw (22) to (23);
						\draw (23) to (32);
						\draw (32) to (37);
						\draw (37) to (36);
						\draw (36) to (35);
						\draw (43) to (46);
						\draw (46) to (47);
						\draw (47) to (48);
						\draw (48) to (57);
						\draw (45) to (52);
						\draw (52) to (53);
						\draw (53) to (54);
						\draw (54) to (55);
						\draw [bend left, looseness=1.25] (3) to (5);
						\draw (57) to (64.center);
						\draw (55) to (63.center);
					\end{pgfonlayer}
				\end{tikzpicture}
			}\caption{\label{fig:bt2}A portion of the graph $G$ corresponding to functions $f$ and $g$ such that $0 \notin \ran(f \upto 4) \cup \ran(g\upto 4)$, $1=f(3)$ and $2=g(2)$.}			
		\end{figure}
		From the given description it is clear that $G$ is $\lDelta{0}{1}$-definable in $f$ and $g$, its degree function $\delta_{G}$ is also $\lDelta{0}{1}$-definable in $f$ and $g$, $\Delta(G)=2$, and $G$ is acyclic. By $\btb 2$, let $c$ be a $2$-coloring of $G$. It is easy to see that if there is $m$ such that $f(m)=n$, then $c(l_{n,0}) \neq c(r_{n,0})$, while if there is $m$ such that $g(m)=n$, then $c(l_{n,0})= c(r_{n,0})$. Thus the set $X=\{n \in \N : c(l_{n,0})=c(r_{n,0})\}$ separates the range of $f$ from the range of $g$.
	\end{proof}
	
	Notice that Lemma \ref{lem:brooksfinitebi} implies that some of the graphs used to prove Proposition \ref{prop:boundedbiwkl} must have at least one infinite connected component. In fact, whenever there is $n$ such that $n \notin \ran(f) \cup \ran(g)$, the graph $G$ built as in the aforementioned proof has at least two infinite connected components (consisting of the vertices $\{l_{n,k}:k \in \N\}$ and $\{r_{n,k}:k \in \N\}$).
	
	\subsection{$\rca$ proves $\btb[3]$}
	
	We now set out to show that $\rca$ can prove $\forall d \geq 3 \,\, \bt[d]$. This is done by clarifying a proof of Tverberg in \cite{Tverberg84}, and carefully setting up the induction therein to obtain a $\rca$ proof.
	
	In this subsection we start by giving a thorough proof of $\btb[3]$ in $\rca$. This is the base case of a slightly subtle induction, which we will discuss in the next subsection. We start with some auxiliary definitions and results which allow for a leaner presentation of the proof of Theorem \ref{thm:bt3}.
	
	\subsubsection{Circle-trees}
	
	\begin{definition}
		Let $G=(V,E)$ and $H=(V',E')$ be graphs such that $V \cap V'= \emptyset$, $\Delta(G)\leq3$ and $\Delta(H) \leq 3$. Let $x \in V$ and $y \in V'$ with $\delta_G(x)=\delta_H(y)=2$. We define the graph $G+_{x,y}H$ as $(V \cup V', E \cup E' \cup \{(x,y)\})$.
	\end{definition}
	
	In what follows, we will only write $G +_{x,y} H$ if the vertex sets of $G$ and $H$ are disjoint and $\delta_G(x)=\delta_H(y)=2$.
	
	\begin{definition}\label{def:cycletrees}
		A finite graph $X$ is a \emph{circle-tree} if there exists $k \in \N$ and a sequence of numbers $(n(i):i<k)$ such that the following conditions hold:
		\begin{itemize}
			\item $n(0)$ codes a circle $C(0)=T(0)$,
			\item for every $0<i<k$, $n(i)$ codes a finite graph $T(i)$ such that $T(i)=T(i-1)+_{x,y}C(i)$ for some circle $C(i)$ with $C(i) \cap T(i-1) = \emptyset$, $x \in T(i-1)$, $y \in C(i)$,
			\item the graph $T(k-1)$ is $X$.
		\end{itemize}
		
		If $k$ and $(n(i):i<k)$ are as above, we say that they \emph{witness} that $T$ is a circle-tree, moreover; we say that the circles $(C(0), \dots, C(k-1))$ \emph{make up the tree $T$ according to $k$ and $(n(i):i<k)$}.
	\end{definition}

	\begin{remark}\label{rem:fuf}
		Giving the right definition of circle-trees for our purposes is quite delicate, as one can easily run into problems with induction.
		Indeed, $\rca$ does not prove that ``finite unions of circles are finite'' (this statement is clearly equivalent to the $\mathsf{FUF}$ principle ``finite unions of finite sets is finite'', which is known to be equivalent to $\mathsf{B}\lSigma{0}{2}$. See \cite[Proposition 6.5.4]{DM}).
		With the definition above, this finiteness is built in. Hence it is not true that every sequence of circles corresponds to a circle-tree.
	\end{remark}
	
	\begin{lemma}\label{lem:subtraction}
		Let $X$ be a circle-tree, as witnessed by $k \geq 2$ and $(n(i):i<k)$, and let $(C(0), \dots, C(k-1))$ be circles which make up $X$ according to $k$ and $(n(i):i<k)$. Then, for every $j \leq k$, $j$ and $(n(i):i<j)$ witness that $T=X \setminus \bigcup_{j \leq i <k }C(i)$ is a circle-tree.
	\end{lemma}
	
	This Lemma follows directly from the definition of circle-trees, and will be instrumental in the rest of our work with these objects.
	
	Note that the formula $\psi(n)$ stating ``$n$ codes a circle-tree'' is $\lDelta{0}{1}$.
	A code for a circle-tree also keeps track of the \emph{construction} of the circle-tree: $k$ is the number of steps in the construction, whereas the numbers $(n(i) : i < k)$ are the instructions to build the circle-tree (which circles to use and how to glue them together).
	Given a finite graph $X$ which happens to be a circle-tree, there may be more than one construction of $X$ as a circle-tree.
	We now show that the construction of a circle-tree is essentially unique. This allows us to refer unambiguously to the elements of a construction of a circle-tree.
	
	To avoid verbosity, we use the following abbreviation: we say that a sequence of circles $(C(i) : i < k)$ \emph{determines} a circle-tree $T$ if there is a sequence $(n(i) : i < k)$ which witnesses that $T$ is a circle-tree and such that $(C(i) : i < k)$ is the sequence which makes up $T$ according to $(n(i) : i < k)$. If $(C(i) : i <k)$ is a sequence of circles which determines the circle-tree $T$, then the set of vertices in $T$ is given by $\bigcup_{i<k}C(i)$. Every circle-tree $T$ is determined by some sequence of circles but specifying a sequence of circles does not fully determine a circle-tree, so we are abusing language here.
	
	\begin{lemma}\label{lem:number}
		Let $T$ be a circle-tree, determined by $(C(i) : i < k)$. Then $|\{x  \in T : \delta(x)=3\}|=2(k-1)$. Consequently, if $T$ is also determined by $(C'(i) : i < k')$, then $k'=k$.
	\end{lemma}
	\begin{proof}
		By an easy induction on $k$.
	\end{proof}
	So we have shown that if $X$ is a circle-tree, then all of its constructions consist of the same number of steps. We set out to show that all of its constructions involve the same set of circles.
	
	\begin{definition}
		Let $T$ be a circle-tree, as witnessed by $k$ and $(n(i) : i < k)$, and let $(C(i): i < k)$ be the circles which make up $T$. The graph $Q_T$ is has vertex set $\{C(i): i < k\}$, and for any $i<j <k$, $C(i)EC(j)$ in $Q_T$ if and only if $C(i) \cup C(j)$ induces a connected subgraph of $T$.
	\end{definition}
	
	\begin{lemma}\label{lem:quotientacyclic}
		Let $T$ be a circle-tree. Then $Q_T$ is acyclic and connected.
	\end{lemma}
	\begin{proof}
		We prove the claim by induction on the number of steps $k$ in a construction of $T$. Formally we are proving $\forall k \,\, \varphi(k)$ where $\varphi(k)$ states ``for every $(n(i) : i < k)$ which codes a circle-tree $T$, the graph $Q_T$ is connected, and for every finite induced subgraph $X \subseteq Q_T$, $X$ is not a cycle''. The formula $\varphi(k)$ is easily seen to be $\lDelta{0}{1}$, so we are allowed to prove $\forall k \,\, \varphi(k)$ by induction in $\rca$.
		
		The base case is $k=1$: in this case $T$ consists of a single circle and $Q_T$ consists of a single node, it is acyclic and connected.
		
		Now assume that $\varphi(k)$ holds and let $T$ be a circle-tree with a construction coded by $k+1$ and $(n(i) : i < k+1)$. Let $(C(i) : i < k+1)$ be the sequence of circles which make up $T$. By Lemma \ref{lem:subtraction} we know that $k$ and $(n(i) : i < k)$ witness that $T'=T \setminus C(k)$ is a circle-tree. We have that $Q_T$ is obtained by $Q_{T'}$ by adding the node $C(k)$ and one edge from $C(k)$ to a previous node.
		This operation preserves connectedness and acyclicity (since $C(k)$ has degree $1$ in $Q_T$).
	\end{proof}
	
	Connected acyclic graphs are known in graph theory as \emph{trees}, and it is well known that they can be constructed starting from an arbitrary node, adding one node at a time in such a way that all partial graphs are connected. This implies that:
	\begin{lemma}\label{lem:arbitraryroot}
		Let $T$ be a circle-tree which is determined by $(C(i) : i < k)$. For any $j < k$ there is some construction $(C'(i) : i < k)$ of $T$ such that $C'(0)=C(j)$.
	\end{lemma}
	
	We now get back to the proof that constructions of circle-trees are essentially unique.
	\begin{lemma}\label{lem:cycleindivisibility}
		Let $X$ be a graph which consists of a single circle and let $C \subseteq X$. Then $C$ is a circle if and only if $C=X$.
	\end{lemma}
	The proof of Lemma \ref{lem:cycleindivisibility} is immediate, considering the degrees of the vertices involved.
	\begin{lemma}\label{lem:acyclicity}
		Let $T$ be a circle-tree which is determined by $(C(i) : i < k)$ and let $C \subseteq T$ be a circle. Then there is $j < k$ such that $C=C(j)$.
	\end{lemma}
	\begin{proof}
		Let $(v(i) : i < n+1)$ be an enumeration of the nodes of $C$ such that $v(i)Ev(i+1)$ for every $i<n$, $|\{v(i) : i < n\}|=n$, and $v(n)=v(0)$. Note that $n > 2$.
		
		If there is some $j$ such that $C \subseteq C(j)$, then by Lemma \ref{lem:cycleindivisibility}, $C=C(j)$ and the claim holds.
		Otherwise, consider the quotient map $\pi \colon T \rightarrow Q_T$. The sequence $(\pi(v(i)) : i < n+1)$ has the property that for every $i<n+1$, either $\pi(v(i))=\pi(v(i+1))$ or $\pi(v(i))E \pi(v(i+1))$. Let $(w(i) : i <m)$ be the subsequence of $(\pi(v(i)) : i < n+1)$ defined as $w(0)=\pi(v(0))$ and $w(i+1)=\pi(v(j))$ where $j$ is least such that $\pi(v(j)) \neq w(i)$. By definition $w(i) \neq w(i+1)$, and moreover, $w(i) \neq w(i+2)$ because otherwise there would be two edges between a vertex in $w(i)$ and a vertex in $w(i+1)$. It is easy to see that $w(0)Ew(m-1)$, so if $(w(i) : i < m)$ is not already a cycle, there must be a minimal $i$ and a minimal $j>i$ such that $w(i)=w(j)$. We have $j \notin \{i+1, i+2\}$, hence $\{w(h) : i \leq h<j\}$ is a cycle in $Q_T$. This contradicts Lemma \ref{lem:quotientacyclic}.
	\end{proof}
	
	\begin{proposition}\label{prop:ctuniqueness}
		Let $T$ be a circle-tree and suppose that $(C(i) : i < k)$ is a sequence of circles which determines $T$. Let $(C'(i) : i < k)$ be a sequence of circles which also determines $T$. Then there is a bijection $f \colon k \rightarrow k$ such that $C(i)=C'(f(i))$ for every $i$.
	\end{proposition}
	\begin{proof}
		By Lemma \ref{lem:acyclicity}, we can define the function $f \colon k \rightarrow k$ as follows: $f(i)$ is defined to be the unique $j$ such that $C'(j)=C(i)$. We claim that $f$ is well defined and is a bijection. To see this, first note that for every $i$, Lemma \ref{lem:acyclicity} guarantees that there is some $j$ such that $C(i)=C'(j)$. Moreover, the map is injective because if $i \neq j$, then $C(i) \cap C(j) = \emptyset$ so by definition $C'(f(i)) \cap C'(f(j)) = \emptyset$. By the finite pigeonhole principle (provable in $\rca$), we conclude that $f$ is also surjective.
	\end{proof}
	
	Proposition \ref{prop:ctuniqueness} guarantees that if $X$ is a finite graph which happens to be a circle-tree, there is a unique finite set of disjoint circles $\{C(i) : i < k\}$ such that the vertex set of $X$ is $\bigcup_{i <k}C(i)$. We will use this fact implicitly throughout the proof of Theorem \ref{thm:bt3}.
	We conclude the section with some other easy results on circle-trees which will prove useful in the proof of Theorem \ref{thm:bt3}.
	
	\begin{lemma}\label{lem:trees}
		If $T$ is a circle-tree, then $3|\{x \in T : \delta(x)=2\}| > |T|$ and hence $2|\{x \in T : \delta(x)=2\}| >|\{x \in T : \delta(x)=3\}|$.
	\end{lemma}
	\begin{proof}
		The first inequality  is obtained by induction on the number of circles $k$ which constitute $T$.
		The second inequality follows from the first one as $T$ only contains nodes of degree two or three.
	\end{proof}
	\begin{lemma}\label{lem:trees2}
		Let $T$ be a circle-tree, $x \in T$. If $\delta(x)=2$, then $T \setminus \{x\}$ is connected, and it is not a circle-tree. If $\delta(x)=3$, then $T \setminus \{x\}$ is not connected, and exactly one of its two connected components is a circle-tree.
	\end{lemma}
	\begin{proof}
		By induction on the number of circles which constitute $T$.
	\end{proof}
	\begin{lemma}\label{lem:treesU}
		Let $T_1$ and $T_2$ be circle-trees, $x \in T_1$ and $y \in T_2$ with $\delta_{T_1}(x)=\delta_{T_2}(y)=2$. Then $T_1+_{x,y}T_2$ is a circle-tree.
	\end{lemma}
	\begin{proof}
		By induction on the number of circles which constitute $T_2$, using Lemma \ref{lem:arbitraryroot}.
	\end{proof}
	\begin{lemma}\label{lem:cyclesintrees}
		Let $G$ be a graph, and let $T \subseteq G$ be a circle-tree such that $|\{t \in T : N(t) \nsubseteq T\}|=1$. Let $x \in T$ and let $C \subseteq G$ be a circle containing $x$. Then $C$ is one of the circles which make up $T$.
	\end{lemma}
	\begin{proof}
		By contradiction assume that $C$ is not one of the circles which make up $T$. By Lemma \ref{lem:acyclicity}, $C$ is not contained in $T$, so $T \cap C$ is a finite graph whose connected components consist of line segments. Hence, there are two nodes $u_0$ and $u_1$ in $T$ such that $\delta_{T\cap C}(u_0)=\delta_{T\cap C}(u_1)=1$. Since there is only one node in $T$ which is connected to $G \setminus T$, it must be that one of the two, say $u_0$, has $\delta_{C}(u_0)=1$. This contradicts the fact that all nodes in a circle have degree $2$.
	\end{proof}
	\subsubsection{$P$-vertices, $Q$-vertices, germs, and circle-trees in graphs with degree $3$}
	
	\begin{definition}
		If $G$ is a graph, we say that $v \in V$ is a \emph{$P$-vertex} of $G$ if $\delta(v)=3$ and there is a circle $C$ in $V$ such that $v \in C$ and $\delta(u) = 2$ for all $u \in C \setminus \{v\}$.
	\end{definition}
	\begin{lemma}\label{lem:pvdisjointcycles}
		Let $G$ be a graph and let $v_0$ and $v_1$ be distinct $P$-vertices of $G$. Let $C_0$ and $C_1$ be circles in $G$ such that $v_i$ is the unique node of degree $3$ in $C_i$. Then $C_0 \cap C_1 = \emptyset$.
	\end{lemma}
	\begin{proof}
		First notice that $C_0 \neq C_1$. Now consider the induced subgraph of $G$ given by $C_0 \cup C_1$, and assume that $C_0 \cap C_1 \neq \emptyset$.
		Hence, $C_0 \cup C_1$ is connected and it is easy to see that there must be some $x \in C_0 \cap C_1$ such that $N_1(x) \nsubseteq C_0 \cap C_1$; hence $\delta(x) \geq 3$. This implies that at least one of $C_0$ and $C_1$ must have two nodes of degree $3$, a contradiction.
	\end{proof}
	
	\begin{lemma}\label{lem:treespv}
		Let $T$ be a circle-tree made up of $k$ circles $(C(i) : i < k)$ according to the sequence $(n(i) : i < k)$. If $k \geq 2$, then there are at least two $P$-vertices in $T$.
		More precisely, there is $i < k-1$ , $u_i \in C_i$ and $u_{k-1} \in C_{k-1}$ such that for $j \in \{i,k-1\}$, $u_j$ is the unique node of degree $3$ within $C_j$.
		
		Moreover, if $k \geq 2$ and $T$ has exactly two $P$-vertices $x$ and $y$, then $Q_T$ is a line segment, and the $P$-vertices of $T$ belong to the circles corresponding to the endpoints of $Q_T$.
	\end{lemma}
	\begin{proof}
		The proof of the first claim hinges on the known graph theoretic fact (which is provable in $\rca$) that a tree with more than one node has at least two nodes of degree $1$. Hence, $Q_T$ has at least two nodes of degree $1$. It is easy to see that the corresponding circles in $T$ contain $P$-vertices.
		
		If $T$ has exactly two $P$-vertices, then $Q_T$ must have only two nodes of degree $1$, and again by basic graph theoretic results this implies that $Q_T$ is a line segment.
	\end{proof}
	
	\begin{definition}
		If $G$ is a graph, we say that $v \in V$ is a \emph{$Q$-vertex} of $G$ if $\delta(v)=3$, $v$ is not a $P$-vertex, but there is a circle $C \subseteq V$ such that $v \in C$, $|\{u \in C: \delta(u) =3\}|\leq 3$ and $N(v) \cap C$ contains at least one vertex of degree $2$.
	\end{definition}
	\begin{remark}
		In $\rca$ one can show, by induction on $k$, that if $G$ is a graph which induces a cycle $(v(i) : i < k)$, then there is a circle in $G$ induced by a subset of the nodes $\{v(i): i < k\}$. This implies that defining $P$-vertices and $Q$-vertices in terms of cycles would yield the same concepts. Indeed this is what Tverberg does in \cite{Tverberg84}\footnote{Tverberg uses the word cycles without an explicit definition. We assume (according to common usage) that his concept agrees with ours.}. In our rendition of his proof we use circles for convenience.
	\end{remark}
	\begin{definition}
		Let $G$ be a graph. We say that a path $(v(i) : i < n)$ is a \emph{germ} in $G$ if $n \geq 3$, $\delta(v(i)) \leq 2$ for all $i< n-1$, and $\delta(v(n-1))=3$.
	\end{definition}
	\begin{lemma}\label{lem:germsintrees}
		Let $G$ be a graph and let $T \subseteq G$ be a circle-tree. Let $P=(v(i) : i < n)$ be a germ in $G$ and assume that $(P \cap T )\setminus \{v(n-1)\} \neq \emptyset$. Then $P \subseteq T$
	\end{lemma}
	\begin{proof}
		If $P \nsubseteq T$ there exists $i<n-1$ such that $v(i) \in T$ and either $v(i-1)$ or $v(i+1)$ does not belong to $T$. Since $\delta_T(v(i)) \geq 2$, we have $\delta_G(v(i)) = 3$, against the definition of germ.
	\end{proof}
	
	\subsubsection{Tverberg's proof}
	
	\begin{theorem}\label{thm:bt3}
		Let $G=(V,E,f)$ be a bounded graph with $\Delta(G) \leq 3$ which does not induce a $4$-clique. Then $G$ is $3$-colorable.
	\end{theorem}
	\begin{proof}
		By Lemma \ref{lem:boundedregular} we can assume without loss of generality that $G$ is $3$-regular.
		
		We assume that $V= \mathbb{N}$ and we describe an effective procedure to define a set $W \subseteq V$ in stages.
		$W$ will be an independent set, and all connected components of $G \setminus W$ will be line segments. Therefore $G \setminus W$ will be $2$-colorable by Lemma \ref{lem:brooksfinitebi}. This shows that $G$ is $3$-colorable.
		
		Assume that $w_0, \dots, w_{i-1} \in V$ have been defined and consider $G_i$ as the subgraph of $G$ induced by vertices $V \setminus \{w_0, \dots, w_{i-1}\}$. Let $s$ be the least vertex of degree $3$ in $G_i$.
		
		\begin{enumerate}[left= 0pt]
			\item \label{it:pv}If $G_i$ has $P$-vertices, let $w_i$ be the first $P$-vertex.
			\item \label{it:germs} If $G_i$ has no $P$-vertices and there is a germ $j_1E j_2 E \dots j_{m-1} E j_m$ with $j_1<s$, we pick the germ with minimal starting node (if there are two germs starting with this node choose the one with minimal second node). Let this germ be $aEb E \dots l E m$.
			\begin{enumerate}
				\item  \label{it:gt}If $G_i \setminus \{m\}$ has a connected component which is a circle-tree, let $T$ be the circle-tree of least code which is a connected component in $G_i \setminus \{m\}$. Choose as $w_i$ the first element of $(N(m) \setminus \{l\}) \cap T$ if this set is nonempty. If $(N(m) \setminus \{l\}) \cap T= \emptyset$, choose $w_i$ as the first element of $(N(m) \setminus \{l\})$.\footnote{We will prove in Claim \ref{claim:notree} that the second case never occurs.}
				\item \label{it:gnt}If $G_i \setminus \{m\}$ does not have a connected component which is a circle-tree, let $w_i=m$.
			\end{enumerate}
			\item \label{it:qv}If $G_i$ has no $P$-vertices nor germs starting with vertices before $s$, if there are $Q$-vertices in $N(s)$, let $w_i$ be the first one. Otherwise, let $w_i=s$.
		\end{enumerate}
		
		\begin{claimn}\label{claim:effectiveness}
			The procedure described above is effective.
		\end{claimn}
		\begin{proof}
			First notice that, for every $i$, and every $v \in G_i$ we can compute (uniformly in $i$) $\delta_{G_i}(v)$ by checking how many of the $G$-neighbors of $v$ belong to $\{w_0, \dots, w_{i-1}\}$. Hence, we can find $s$ effectively. Moreover, by $\bDelta{0}{1}$-induction on $i$, we have $\{v \in G_i : \delta_{G_i}(v) < 3\} \subseteq \bigcup_{k < i}N(w_k)$ and in particular $|\{v \in G_i : \delta_{G_i}(v) < 3\}| \leq 3i$.
			
			The algorithm first needs to check whether $G_i$ has $P$-vertices (cf.\ case \ref{it:pv}). For that, it suffices to look at all vertices of $G_i$ of degree $2$ and their neighbors, and check whether a subset of these vertices forms a circle which only contains one vertex of degree $3$. This procedure is effective as $\{v \in G_i : \delta_{G_i}(v) < 3\}$ is finite.
			
			Checking whether or not $G_i$ has germs which begin with nodes less than $s$ (cf.\ case \ref{it:germs}) is also computable: similarly to the case of finding $P$-vertices, germs consist of several nodes of degree less than $3$. Since there are only finitely many such nodes, we can exhaustively check all possibilities.
			
			If we find a germ $a E b E \dots l E m$ with $a < s$ minimal, our next task is to check whether $G_i \setminus \{m\}$ has a connected component which is a circle-tree.
			This is possible due to Lemma \ref{lem:trees}, which guarantees that in any circle-tree, at least a third of the nodes have degree $2$.
			To see why this is the case, let $H=G_i \setminus \{m\}$. For every $x \in H$ with $\delta_{H}(x)=2$ we consider the set $N_n(x)$ for increasingly large values of $n$, until we find either
			\begin{itemize}
				\item some $n$ such that $N_n(x)=N_{n+1}(x)$, or
				\item some $n$ such that $|N_n(x)| > 3|\{x \in H : \delta_H(x)=2 \}|$.
			\end{itemize}
			It is clear that one of the two alternatives must eventually come about. In the first case, we have that $N_n(x)$ is a connected component of $H$, and we can simply check whether or not it is a circle-tree. In the second case the connected component of $x$ cannot be a circle-tree.
			
			If $G_i$ has no $P$-vertices nor germs starting with vertices less than $s$, we are in case \ref{it:qv}, so we check whether there are $Q$-vertices in $N(s)$. This check is again effective, because a $Q$ vertex must lie on a circle with less than $3$ nodes of degree $3$, and there are only finitely many nodes of degree $2$.
		\end{proof}
		By Claim \ref{claim:effectiveness} the function $g : i \mapsto w_i$ which gives us the points we pick out of $G$ with this construction is $\lDelta{0}{1}$-definable from $G$.
		By the end of this proof we will see that the range of the sequence $(w_i : i \in \N)$ exists as a set, even in $\rca$.
		Indeed, we will prove that $W$ is independent and $G \setminus W$ has degree $2$, so that $v \in W$ if and only if $v$ is enumerated by the function $g$, if and only if no neighbor of $v$ in $G$ is ever enumerated by $g$. This gives a $\lDelta{0}{1}(G)$ definition of $W$. In the meantime, we will refer to $W=\{w_i: i \in \N\}$ and $G \setminus W$ informally (as abbreviations) to give a better intuition of how the proof works.
		
		It remains to show that no connected component of $G \setminus W$ is a circle, $W$ is independent, $\Delta(G \setminus W)=2$, and all connected components of $G \setminus W$ are finite (hence, the connected components are line segments). We will prove each of these facts via a sequence of claims. Claim \ref{claim:notree} states that, at each stage $i$, no connected component of the graph $G_i$ is a circle-tree. It is immediate to see that this implies that no connected component of $G \setminus W$ is a circle. Claim \ref{claim:independent} readily implies that $G \setminus W$ is independent. Claim \ref{claim:germs} is a technical lemma on germs which shows that the construction works as intended. We note that the precise computational bound we obtain there is actually necessary to make sure the inductive argument is formalizable in $\rca$. Claims \ref{claim:rightdegree} and \ref{claimcompo} state that $\Delta(G \setminus W)=2$ and that the connected components of $G \setminus W$ are line segments.
		
		\begin{claimn}\label{claim:notree}
			For every $i$, no connected component of $G_i$ is a circle-tree. Therefore, no component of $G \setminus W$ is a circle-tree.
		\end{claimn}
		\begin{proof}
			Since, given a bounded graph $H$ and a finite $X \subseteq H$, the statement ``$X$ is circle-tree'' is $\lDelta{0}{1}(H)$, the statement ``$H$ has a component which is a circle-tree'' is $\lSigma{0}{1}(H)$, being equivalent to:
			\[
			\exists x \exists n (N_n(x)=N_{n+1}(x) \land N_n(x) \text{ is a circle-tree}).
			\]
			
			Therefore, we can prove that $G_i$ does not have a component which is a circle-tree by $\bPi{0}{1}$-induction on $i$.
			
			For the base case, $G=G_0$ is $3$-regular, while circle-trees have nodes of degree $2$, hence no connected component of $G$ is a circle-tree.
			
			Now assume $G_i$ does not have any connected component which is a circle-tree and recall $G_{i+1}=G_i \setminus \{w_i\}$. Suppose by contradiction that there is a component, call it $T$, of $G_{i+1}$ which is a circle-tree. Inductively $T$ was not a component of $G_i$, hence $w_i \in N_1(T) \setminus T$. Now there are three cases:
			\begin{enumerate}[left= 0pt]
				\item In case \ref{it:pv}, $w_i$ is a $P$-vertex of $G_i$.
				Let $y \in T \cap N(w_i)$ (notice that $\delta_H(y)=2$). Since $w_i$ is a $P$-vertex in $G_i$, there must be a circle $C \subseteq G_i$ such that every node in $C \setminus \{w_i\}$ has degree $2$ in $G_i$. Therefore, we can see that $C \setminus \{w_i\}$ is a connected component of $G_{i+1}$, and in particular $C \cap T = \emptyset$. Putting everything together we get that $T +_{y,w_i} C$ is a connected component of $G_i$ which, by Lemma \ref{lem:treesU}, is a circle-tree, against our inductive assumption.
				\item If $w_i$ is picked following case \ref{it:germs}, there are two possibilities:
				\begin{enumerate}
					\item \label{it:alwaysinT} In case \ref{it:gt}, there is a node $m \in G_i$ such that $G_i \setminus \{m\}$ has at least one connected component which is a circle-tree, and $w_i \in N(m)$. So, let $T' \subseteq G_i \setminus \{m\}$ be the circle-tree of least code which is a connected component.
					
					First, notice that a priori it could be that $T' \cap (N(m) \setminus \{l\})= \emptyset$, and hence $w_i$ does not belong to $T'$.
					We show that this can never actually happen: indeed, by the inductive hypothesis, $T'$ is not a connected component of $G_i$, so $m \in N_1(T') \setminus T'$. Since, $\delta_{G_i \setminus \{m\}}(l)=1$, and circle-trees do not have nodes of degree $1$, $l \notin T'$ and hence $(N(m) \setminus \{l\}) \cap T' \neq \emptyset$. So $w_i \in T'$, and we get
					$3 \geq \delta_{G_i}(w_i)>\delta_{T'}(w_i) \geq 2$, so $\delta_{T'}(w_i)=2$ and by Lemma \ref{lem:trees2}, $T' \setminus \{w_i\}=T' \cap G_{i+1}$ is connected and it is not a circle-tree. More precisely, if $T'$ is determined by circles $(C(i) : i < k)$, then there is some $j <k$ such that $w_i \in C(j)$, and hence $C(j) \cap G_{i+1}$ is a line segment.
					
					Since $T$ is a connected component of $G_{i+1}$ and $T'$ is connected, either $T' \subseteq T$ or $T \cap T' = \emptyset$. Since $T'$ has nodes of degree $1$ (namely, the neighbors of $w_i$ in $C(j)$), we have $T' \nsubseteq T$, so $T \cap T' = \emptyset$.
					
					Since $N(w_i) \cap T \neq \emptyset$, $w_i \in T'$ and $w_i$ has at least two neighbors in $T'$, it follows that $\{m\}=N(w_i) \cap T$, and hence $\delta_T(m)=2$. But then, by Lemma \ref{lem:treesU}, $T+_{w_i,v_m}T'$ is a circle-tree that is a component of $G_i$. This goes against the inductive assumption.
					\item if $w_i$ was picked following rule \ref{it:gnt}, it follows from the definition that $G_{i+1}$ does not have any connected components which are circle-trees,
				\end{enumerate}
				\item lastly, if $w_i$ is picked applying rule \ref{it:qv}, then $w_i=s$ if there are no $Q$-vertices in $N(s)$, otherwise $w_i$ is the least $Q$-vertex in $N(s)$.
				
				In case $w_i$ is a $Q$-vertex, by definition, $w_i$ belongs to some circle $C$ and $N(w_i) \cap C$ contains at least one node $x$ with $\delta_{G_i}(x)=2$. Since $\delta_{G_{i+1}}(x)=1$, $x \notin T$, and since $C \setminus \{w_i\}$ is connected, $C \cap T = \emptyset$. Again by the inductive assumption $N(w_i) \cap T \neq \emptyset$ so we have $|N(w_i) \cap T|=1$.
				
				If $T$ consists of only one circle, and $u \in N(w_i) \cap T$, then $u$ is a $P$-vertex in $G_i$. This means that we would be in case \ref{it:pv} and pick $w_i=u$, rather than follow rule \ref{it:qv}. This is a contradiction.
				
				If $T$ consists of more than one circle, by Lemma \ref{lem:treespv}, $T$ has at least two $P$-vertices, say $u_0$ and $u_1$. Since $u_0$ and $u_1$ are $P$-vertices, there are disjoint circles $C_0 \subseteq T$ and $C_1 \subseteq T$ such that $u_0$ is the only node of $C_0$ of degree $3$ and $u_1$ is the only node of $C_1$ of degree $3$.  Since $|N(w_i) \cap T|=1$, it follows that there is $j \in \{0,1\}$ such that $N(w_i) \cap C_j= \emptyset$. Hence $u_j$ is a $P$-vertex in $G_i$, and again we reach a contradiction.
				
				Lastly, assume $w_i$ is not a $Q$-vertex of $G_i$. If $T$ consists of one circle, then, if $|N(w_i) \cap T|=1$, we would have had a $P$-vertex in $G_i$, which is a contradiction.
				If $1 < |N(w_i) \cap T| \leq 3$, then at least two of the the neighbors of $w_i$ in $T$ would be $Q$-vertices in $G_i$, which is again a contradiction.
				
				If $T$ consists of more than one circle, then $T$ has at least two $P$-vertices $u_0$ and $u_1$. Let $C_0$ and $C_1$ be circles to which $u_0$ and $u_1$, respectively, belong, and such that $u_j$ is the only vertex of $C_j$ with $\delta_T(u_j)=3$. Since there are no $P$-vertices in $G_i$, it follows that $N(w_i) \cap C_j \neq \emptyset$ for $j \in \{0,1\}$. Let $t_0$ and $t_1$ be elements of $N(w_i) \cap C_0$ and $N(w_i) \cap C_1$, respectively. For $j \in \{0,1\}$ we have: $\delta_{G_i}(t_j)=3$, $|\{x \in C_j : \delta_{G_i}(x)=3\}| \leq 3$, and clearly at least one of the $t_j$ has a neighbor in $C_j$ of degree $2$. Hence each $t_j$ is a $Q$-vertex in $G_i$. This is again a contradiction.\qedhere
			\end{enumerate}
		\end{proof}
		\begin{claimn}\label{claim:independent}
			For all $i$, $w_i$ has degree $3$ in $G_i$. In particular $w_i \notin N(w_j)$ for all $j<i$, and hence $W$ is independent.
		\end{claimn}
		
		\begin{proof}
			If we are in case \ref{it:pv}, \ref{it:gnt}, or \ref{it:qv} this is immediate.
			
			If we choose $w_i$ according to rule \ref{it:gt}, we found a germ $aE \dots ElEm$ in $G_i$ such that a component of $G_i \setminus \{m\}$ is a circle-tree. In this case, let $T \subseteq G_i \setminus \{m\}$ be the circle-tree with least code. By the proof of point \ref{it:alwaysinT} in Claim \ref{claim:notree}, $w_i \in T \cap N(m)$. Since any node in a circle-tree has degree at least $2$ and $w_i \in N(m)$, we have that $w_i$ has degree $3$ in $G_i$ also in this case.
		\end{proof}
		
		\begin{claimn}\label{claim:germs}
			For every $i$, if we are in case \ref{it:germs} and $aEbE\dots Em$ is the germ in $G_i$ which has the least starting node (with $b$ least if there are two germs in $G_i$ starting with $a$), then there is no germ which starts with $aEb$ in $G_{4i}$. Indeed the connected component of $b$ in $G_{4i} \setminus \{a\}$ is a line segment.
		\end{claimn}
		\begin{proof}
			First, we clarify formally what we want to prove.
			
			Let $\psi(i)$ be the formula ``there are no $P$-vertices in $G_i$'', let $\varphi(a,b,n,i)$ be the formula ``$n$ codes a germ in $G_i$ which starts with $aEb$'', and let $\rho(a,b,n,i)$ be the formula ``$n$ codes a germ in $G_i$ which starts with $aEb$, $a$ is the minimal node which starts a germ in $G_i$ and, if $a$ is the starting node of two germs in $G_i$, then $b$ is the least neighbor of $a$''. The first part of the statement is
			\[\forall i \forall a \forall b[(\psi(i) \land \exists n \rho(a,b,n,i)) \rightarrow \forall m \neg \varphi(a,b,m,4i)]\footnote{By considerations similar to those in the proof of Claim \ref{claim:effectiveness}, all quantifiers in the formula $\chi(i)=\forall a \forall b[(\psi(i) \land \exists n \rho(a,b,n,i)) \rightarrow \forall m \neg \varphi(a,b,m,4i)$ are actually bounded, so we can prove $\forall i \,\, \chi(i)$ by $\bDelta{0}{1}$-induction.}.\]
			
			Our proof actually shows that, if $\psi(i)$ holds and there are $a$, $b$, and $n$ such that $\rho(a,b,n,i)$, then the connected component of $b$ in $G_{4i} \setminus \{a\}$ is a line segment.
			
			Since $G$ is $3$-regular, it has no germs, so the claim holds for $i=0$.
			
			So, suppose that $i \geq 1$, $G_i$ has no $P$-vertices,  $s$ is the least node with $\delta_{G_i}(s)=3$, $a$ is the minimal node below $s$ which starts a germ in $G_i$, and, if $a$ starts more than one germ in $G_i$, then $b$ is also least among the neighbors of $a$. Let $aEbE \dots ElEm$ be the unique germ in $G_i$ starting with $aEb$.
			
			First note that, if $G_i \setminus \{m\}$ does not contain a connected component which is a circle-tree, then $w_i=m$, and hence $G_{i+1}$ does not contain a germ starting with $aEb$. Moreover, the unique maximal path in $G_{i+1}$ starting with $aEb$ is given by $aEbE\dots El$ and it has a terminal node, namely $l$, with degree $1$. Claim \ref{claim:independent} gives that, also in $G_{4i}$, there are no germs starting with $aEb$ and the line segment $P$ given by $bE \dots El$ is the connected component of $b$ in $G_{4i} \setminus \{a\}$.
			
			Suppose now that $G_i \setminus \{m\}$ does contain a connected component which is a circle-tree, and let $T$ be one such component.
			We already know from the proof of Claim \ref{claim:notree} that $w_i \in T \cap N(m)$. Reasoning on the $P$-vertices of $T$ (cf.\ Lemma \ref{lem:treespv}), since there are no $P$-vertices in $G_i$, it is easy to see that $\{x,y\}=N(m) \setminus \{l\} \subseteq T$ and, if $T$ consists of more than one circle, then $Q_T$ is a line segment and $x$ and $y$ belong to the circles corresponding to the endpoints of $Q_T$. Let $x < y$.
			
			First, assume that $T$ consists of only one circle (case A). There are two subcases:
			\begin{enumerate}[(i),left = 0pt]
				\item if $xEy$, then the situation is as depicted in Figure \ref{fig:onecycleconn}.
				\begin{figure}
					\centering
					\begin{tikzpicture}
						\begin{pgfonlayer}{nodelayer}
							\node [style=empty dot] (0) at (-9, 0) {$a$};
							\node [style=empty dot] (1) at (-8, -1) {};
							\node [style=empty dot] (2) at (-8, 1) {};
							\node [style=none] (3) at (-7, 1) {};
							\node [style=empty dot] (4) at (-4, 1) {$m$};
							\node [style=empty dot] (5) at (-3, 2) {$x$};
							\node [style=empty dot] (7) at (-1, 2) {$z$};
							\node [style=empty dot] (8) at (-3, 0) {$y$};
							\node [style=empty dot] (9) at (-5, 1) {$l$};
							\node [style=none] (12) at (-6, 1) {};
							\node [style=none] (13) at (-6.5, 1) {$\dots$};
							\node [style=empty dot] (14) at (-1, 0) {};
							\node [style=empty dot] (15) at (-9, -4) {$a$};
							\node [style=empty dot] (16) at (-8, -5) {};
							\node [style=empty dot] (17) at (-8, -3) {};
							\node [style=none] (18) at (-7, -3) {};
							\node [style=empty dot] (19) at (-4, -3) {$m$};
							\node [style=empty dot] (21) at (-1, -2) {$z$};
							\node [style=empty dot] (22) at (-3, -4) {$y$};
							\node [style=empty dot] (23) at (-5, -3) {$l$};
							\node [style=none] (24) at (-6, -3) {};
							\node [style=none] (25) at (-6.5, -3) {$\dots$};
							\node [style=empty dot] (26) at (-1, -4) {};
							\node [style=none] (27) at (-11, 0) {$G_i$};
							\node [style=none] (28) at (-11, -4) {$G_{i+1}$};
						\end{pgfonlayer}
						\begin{pgfonlayer}{edgelayer}
							\draw (1) to (0);
							\draw (0) to (2);
							\draw (2) to (3.center);
							\draw (9) to (4);
							\draw (4) to (8);
							\draw (5) to (7);
							\draw (4) to (5);
							\draw (9) to (12.center);
							\draw (5) to (8);
							\draw (8) to (14);
							\draw (14) to (7);
							\draw (16) to (15);
							\draw (15) to (17);
							\draw (17) to (18.center);
							\draw (23) to (19);
							\draw (19) to (22);
							\draw (23) to (24.center);
							\draw (22) to (26);
							\draw (26) to (21);
						\end{pgfonlayer}
					\end{tikzpicture}
					\caption{\label{fig:onecycleconn}
						An application of the rule \ref{it:gt}, in  subcase A(i).}
				\end{figure}
				We have that $\delta_{G_{i+1}}(y)=2$, and the unique path starting with $aEb \dots$ has a terminal node with degree $1$, namely $z$, the unique element in $N(x) \setminus \{m,y\}$. Hence there is no germ starting with $aEb$ in $G_{i+1}$ and the connected component of $b$ in $G_{i+1} \setminus \{a\}$ is a line segment. Again since $G_{4i} \subseteq G_{i+1}$, there is also no germ in $G_{4i}$ stating with $aEb$ and the connected component of $b$ in $G_{4i} \setminus \{a\}$ is a line segment.
				\item if $x$ and $y$ are not adjacent, then the situation is as depicted in Figure \ref{fig:onecycledisc}.
				\begin{figure}
					\centering
					\begin{tikzpicture}
						\begin{pgfonlayer}{nodelayer}
							\node [style=empty dot] (0) at (-9, 0) {$a$};
							\node [style=empty dot] (1) at (-8, -1) {};
							\node [style=empty dot] (2) at (-8, 1) {$b$};
							\node [style=none] (3) at (-7, 1) {};
							\node [style=empty dot] (4) at (-4, 1) {$m$};
							\node [style=empty dot] (5) at (-2, 2) {$x$};
							\node [style=empty dot] (7) at (-1, 1.5) {$z$};
							\node [style=empty dot] (8) at (-2, 0) {$y$};
							\node [style=empty dot] (9) at (-5, 1) {$l$};
							\node [style=none] (12) at (-6, 1) {};
							\node [style=none] (13) at (-6.5, 1) {$\dots$};
							\node [style=empty dot] (14) at (-1, 0.5) {$g$};
							\node [style=empty dot] (15) at (-2.5, 1) {$f$};
							\node [style=empty dot] (16) at (-9, -4) {$a$};
							\node [style=empty dot] (17) at (-8, -5) {};
							\node [style=empty dot] (18) at (-8, -3) {$b$};
							\node [style=none] (19) at (-7, -3) {};
							\node [style=empty dot] (20) at (-4, -3) {$m$};
							\node [style=empty dot] (22) at (-1, -2.5) {$z$};
							\node [style=empty dot] (23) at (-2, -4) {$y$};
							\node [style=empty dot] (24) at (-5, -3) {$l$};
							\node [style=none] (25) at (-6, -3) {};
							\node [style=none] (26) at (-6.5, -3) {$\dots$};
							\node [style=empty dot] (27) at (-1, -3.5) {$g$};
							\node [style=empty dot] (28) at (-2.5, -3) {$f$};
							\node [style=empty dot] (29) at (-9, -8) {$a$};
							\node [style=empty dot] (30) at (-8, -9) {};
							\node [style=empty dot] (31) at (-8, -7) {$b$};
							\node [style=none] (32) at (-7, -7) {};
							\node [style=empty dot] (33) at (-4, -7) {$m$};
							\node [style=empty dot] (34) at (-1, -6.5) {$z$};
							\node [style=empty dot] (36) at (-5, -7) {$l$};
							\node [style=none] (37) at (-6, -7) {};
							\node [style=none] (38) at (-6.5, -7) {$\dots$};
							\node [style=empty dot] (39) at (-1, -7.5) {$g$};
							\node [style=empty dot] (40) at (-2.5, -7) {$f$};
							\node [style=none] (41) at (-11, 0) {$G_i$};
							\node [style=none] (42) at (-11, -4) {$G_{i+1}$};
							\node [style=none] (43) at (-11, -8) {$G_{i+2}$};
						\end{pgfonlayer}
						\begin{pgfonlayer}{edgelayer}
							\draw (1) to (0);
							\draw (0) to (2);
							\draw (2) to (3.center);
							\draw (9) to (4);
							\draw (4) to (8);
							\draw (5) to (7);
							\draw (4) to (5);
							\draw (9) to (12.center);
							\draw (8) to (14);
							\draw (14) to (7);
							\draw (8) to (15);
							\draw (15) to (5);
							\draw (17) to (16);
							\draw (16) to (18);
							\draw (18) to (19.center);
							\draw (24) to (20);
							\draw (20) to (23);
							\draw (24) to (25.center);
							\draw (23) to (27);
							\draw (27) to (22);
							\draw (23) to (28);
							\draw (30) to (29);
							\draw (29) to (31);
							\draw (31) to (32.center);
							\draw (36) to (33);
							\draw (36) to (37.center);
							\draw (39) to (34);
						\end{pgfonlayer}
					\end{tikzpicture}
					\caption{\label{fig:onecycledisc} An application of the rule \ref{it:gt}, in subcase A(ii).}
				\end{figure}
				Since $G_{i+1}=G_i \setminus \{x\}$, there still are no $P$-vertices in $G_{i+1}$ and, once again, the least germ in $G_{i+1}$ is $aEbE \dots EmEy$. Let $\{f,g\}=N(y) \setminus \{m\}$.
				
				We claim that no connected component of $G_{i+1} \setminus \{y\}$ is a circle-tree. Indeed, if this were the case then $N(y) \cap U \neq \emptyset$, but $U$ cannot contain $m$, because $\delta_{G_{i+1} \setminus \{y\}}(m)=1$ and similarly, the connected components of both $f$ and $g$ contain some node of degree $1$ (those nodes which were adjacent to $x$ in $G_i$).
				
				Therefore, in $G_{i+1}$ we apply rule \ref{it:gnt}, deleting $y$. In $G_{i+2}$, there are no germs starting with $aEb$ and the unique path starting with $aEb$ has a terminal node, namely $m$, with degree $1$ in $G_{i+2}$. Again we have that $aEb$ does not start a germ in $G_{4i}$ and the connected component of $b$ in $G_{4i} \setminus \{a\}$ is a line segment.
			\end{enumerate}

			\begin{figure}
				\centering
				\begin{tikzpicture}
					\begin{pgfonlayer}{nodelayer}
						\node [style=empty dot] (0) at (-9, 0) {$a$};
						\node [style=empty dot] (1) at (-8, -1) {};
						\node [style=empty dot] (2) at (-8, 1) {};
						\node [style=none] (3) at (-7, 1) {};
						\node [style=empty dot] (4) at (-4, 1) {$m$};
						\node [style=empty dot] (5) at (-3, 0) {$y$};
						\node [style=empty dot] (8) at (-3, 2) {$x$};
						\node [style=empty dot] (9) at (-5, 1) {$l$};
						\node [style=none] (12) at (-6, 1) {};
						\node [style=none] (13) at (-6.5, 1) {$\dots$};
						\node [style=empty dot] (14) at (-2, 3) {};
						\node [style=empty dot] (15) at (-2, 2) {};
						\node [style=empty dot] (16) at (-1.5, 1.5) {};
						\node [style=empty dot] (17) at (-3, 3) {};
						\node [style=empty dot] (18) at (-0.5, 1.5) {};
						\node [style=empty dot] (19) at (-1.5, 0.5) {};
						\node [style=empty dot] (20) at (-0.5, 0.5) {};
						\node [style=empty dot] (21) at (-2, 0) {};
						\node [style=empty dot] (22) at (-2, -1) {};
						\node [style=empty dot] (23) at (-9, -5) {$a$};
						\node [style=empty dot] (24) at (-8, -6) {};
						\node [style=empty dot] (25) at (-8, -4) {};
						\node [style=none] (26) at (-7, -4) {};
						\node [style=empty dot] (27) at (-4, -4) {$m$};
						\node [style=empty dot] (28) at (-3, -5) {$y$};
						\node [style=empty dot] (30) at (-5, -4) {$l$};
						\node [style=none] (31) at (-6, -4) {};
						\node [style=none] (32) at (-6.5, -4) {$\dots$};
						\node [style=empty dot] (33) at (-2, -2) {};
						\node [style=empty dot] (34) at (-2, -3) {};
						\node [style=empty dot] (35) at (-1.5, -3.5) {};
						\node [style=empty dot] (36) at (-3, -2) {};
						\node [style=empty dot] (37) at (-0.5, -3.5) {};
						\node [style=empty dot] (38) at (-1.5, -4.5) {};
						\node [style=empty dot] (39) at (-0.5, -4.5) {};
						\node [style=empty dot] (40) at (-2, -5) {};
						\node [style=empty dot] (41) at (-2, -6) {};
						\node [style=none] (42) at (-11, 0) {$G_i$};
						\node [style=none] (43) at (-11, -5) {$G_{i+1}$};
					\end{pgfonlayer}
					\begin{pgfonlayer}{edgelayer}
						\draw (1) to (0);
						\draw (0) to (2);
						\draw (2) to (3.center);
						\draw (9) to (4);
						\draw (4) to (8);
						\draw (4) to (5);
						\draw (9) to (12.center);
						\draw (8) to (17);
						\draw (17) to (14);
						\draw (14) to (15);
						\draw (8) to (15);
						\draw (15) to (16);
						\draw (16) to (19);
						\draw (16) to (18);
						\draw (18) to (20);
						\draw (20) to (19);
						\draw (5) to (21);
						\draw (5) to (22);
						\draw (22) to (21);
						\draw (21) to (19);
						\draw (24) to (23);
						\draw (23) to (25);
						\draw (25) to (26.center);
						\draw (30) to (27);
						\draw (27) to (28);
						\draw (30) to (31.center);
						\draw (36) to (33);
						\draw (33) to (34);
						\draw (34) to (35);
						\draw (35) to (38);
						\draw (35) to (37);
						\draw (37) to (39);
						\draw (39) to (38);
						\draw (28) to (40);
						\draw (28) to (41);
						\draw (41) to (40);
						\draw (40) to (38);
					\end{pgfonlayer}
				\end{tikzpicture}
				\caption{\label{fig:morecycles} The first step of an application of the rule \ref{it:gt}, in case B.}
			\end{figure}
			
			The second option (case B) is that $T$ consists of more than one circle. As mentioned above, the situation must be as in Figure \ref{fig:morecycles} (i.e.\ $Q_T$ must be a line segment and $m$ must have one neighbor in each of the circles corresponding to the endpoints of $Q_T$.).
			We apply rule \ref{it:gt}, so we pick $w_i=x$ and define $G_{i+1}=G_i \setminus \{x\}$. We claim that there are no $P$-vertices in $G_{i+1}$.
			
			Towards a contradiction, suppose that there is a $P$-vertex in $G_{i+1}$. Let $t$ be the least one. By definition, there must be a circle $C \subseteq G_{i+1}$ such that $t \in C$, $\delta_{G_{i+1}}(t)=3$, and $t$ is the only node of $C$ which has degree $3$ in $G_{i+1}$. Since there were no $P$-vertices in $G_i$, there must be at least one $r \in C \setminus \{t\}$ with $\delta_{G_{i}}(r)=3$.
			
			Since $\delta_{G_{i+1}}(r)=2$, it follows that $rEx$, so there are two possibilities: either $r \in T$ or $r=m$. Since $mEy$, and $m$ has degree $2$ in $G_{i+1}$, it follows that, if $r=m$, then $y \in C$. So in both cases, we have that $C \cap T \neq \emptyset$.
			
			Now consider the graph $G'$ obtained by deleting the edge $(m,x)$ from $G_i$. Clearly $C$ is a circle in $G'$ too. In $G'$, we can appeal to Lemma \ref{lem:cyclesintrees} to conclude that $C$ is one of the circles which make up $T$. Now it is easy to see that, among the circles which make up $T$, the circle containing $x$ is no longer a circle in $G_{i+1}$, and all other circles have two nodes of degree $3$. So $t$ cannot be a $P$-vertex in $G_{i+1}$.
			
			Given that $G_{i+1}$ contains the germ $aEb \dots EmEy$, and that $a$ is less than the first element of $G_{i+1}$ of degree $3$, we pick $w_{i+1}$ following rule \ref{it:germs}.
			Now suppose that $a'Eb'E \dots El'Em'$ is the germ in $G_{i+1}$ which is taken care of at stage $i+1$. If $m'=y$ (case B(i)), then $w_{i+1}=y$ and we are done. This is shown in Figure \ref{fig:endimmediately}.
			
			\begin{figure}
				\centering
				\begin{tikzpicture}
					\begin{pgfonlayer}{nodelayer}
						\node [style=empty dot] (0) at (-9, 0) {$a$};
						\node [style=empty dot] (1) at (-8, -1) {};
						\node [style=empty dot] (2) at (-8, 1) {};
						\node [style=none] (3) at (-7, 1) {};
						\node [style=empty dot] (4) at (-4, 1) {$m$};
						\node [style=empty dot] (5) at (-3, 0) {$y$};
						\node [style=empty dot] (9) at (-5, 1) {$l$};
						\node [style=none] (12) at (-6, 1) {};
						\node [style=none] (13) at (-6.5, 1) {$\dots$};
						\node [style=empty dot] (14) at (-2, 3) {};
						\node [style=empty dot] (15) at (-2, 2) {};
						\node [style=empty dot] (16) at (-1.5, 1.5) {};
						\node [style=empty dot] (17) at (-3, 3) {};
						\node [style=empty dot] (18) at (-0.5, 1.5) {};
						\node [style=empty dot] (19) at (-1.5, 0.5) {};
						\node [style=empty dot] (20) at (-0.5, 0.5) {};
						\node [style=empty dot] (21) at (-2, 0) {};
						\node [style=empty dot] (22) at (-2, -1) {};
						\node [style=empty dot] (23) at (-9, -5) {$a$};
						\node [style=empty dot] (24) at (-8, -6) {};
						\node [style=empty dot] (25) at (-8, -4) {};
						\node [style=none] (26) at (-7, -4) {};
						\node [style=empty dot] (27) at (-4, -4) {$m$};
						\node [style=empty dot] (30) at (-5, -4) {$l$};
						\node [style=none] (31) at (-6, -4) {};
						\node [style=none] (32) at (-6.5, -4) {$\dots$};
						\node [style=empty dot] (33) at (-2, -2) {};
						\node [style=empty dot] (34) at (-2, -3) {};
						\node [style=empty dot] (35) at (-1.5, -3.5) {};
						\node [style=empty dot] (36) at (-3, -2) {};
						\node [style=empty dot] (37) at (-0.5, -3.5) {};
						\node [style=empty dot] (38) at (-1.5, -4.5) {};
						\node [style=empty dot] (39) at (-0.5, -4.5) {};
						\node [style=empty dot] (40) at (-2, -5) {};
						\node [style=empty dot] (41) at (-2, -6) {};
						\node [style=none] (42) at (-11, 0) {$G_{i+1}$};
						\node [style=none] (43) at (-11, -5) {$G_{i+2}$};
					\end{pgfonlayer}
					\begin{pgfonlayer}{edgelayer}
						\draw (1) to (0);
						\draw (0) to (2);
						\draw (2) to (3.center);
						\draw (9) to (4);
						\draw (4) to (5);
						\draw (9) to (12.center);
						\draw (17) to (14);
						\draw (14) to (15);
						\draw (15) to (16);
						\draw (16) to (19);
						\draw (16) to (18);
						\draw (18) to (20);
						\draw (20) to (19);
						\draw (5) to (21);
						\draw (5) to (22);
						\draw (22) to (21);
						\draw (21) to (19);
						\draw (24) to (23);
						\draw (23) to (25);
						\draw (25) to (26.center);
						\draw (30) to (27);
						\draw (30) to (31.center);
						\draw (36) to (33);
						\draw (33) to (34);
						\draw (34) to (35);
						\draw (35) to (38);
						\draw (35) to (37);
						\draw (37) to (39);
						\draw (39) to (38);
						\draw (41) to (40);
						\draw (40) to (38);
					\end{pgfonlayer}
				\end{tikzpicture}
				\caption{\label{fig:endimmediately} The second step of an application of rule \ref{it:gt}, in case B(i).}
			\end{figure}
			
			If instead $m' \neq y$ (case B(ii)),	we have $a'< a$.
			Since $a$ is the minimal node which starts a germ in $G_i$, $a'Eb'E \dots El'Em'$ is not a germ in $G_{i}$. Therefore, it must be that some nodes among $\{a', \dots, l'\}$ have gone from having degree $3$ to having degree $2$ going from $G_i$ to $G_{i+1}$. So, let $q \in \{a', \dots, l'\} \cap N(w_i)$. Since $w_i=x$, we have that either $q=m$, or $q \in T$.
			
			We claim that if $q = m$, then $m = a'$, and that this leads to a contradiction. Indeed, if $P \subseteq G_{i+1}$ is a germ with $m \in P$ which does not have $y$ as a terminal node, then, since $\delta_{G_{i+1}}(y)=3$, we have $y \notin P$. So the only option is that $m$ is the first node of $P$, so $a'=m$. Now by assumption, if $s$ denotes the least element of $G_{i}$ with $\delta_{G_i}(s)=3$, we have $a < s$. But $m$ has degree $3$ in $G_i$, hence $a< s \leq m$, contradicting $a'<a$.
			
			So, the only remaining option is that $q \in T \cap N(x)$. Since $q$ has degree $3$ in $G_i$, it must be that $q$ is the unique node of $C(0)$ which is connected to some node in $C(1)$. Since $q$ is not the terminal node of the germ under consideration, Lemma \ref{lem:germsintrees} implies that $\{a',b' \dots, l',m'\} \subseteq T$. In particular, we have $q=l'$, $\{a', b' \dots, l'\} \subseteq C(0)$, and $m' \in C(1)$. Applying rule \ref{it:germs} to this germ leads to a situation in which the line segment $C(0) \cap G_{i+2}=C(0) \setminus \{x\}$ is a connected component in $G_{i+2}$. This case is shown in Figure \ref{fig:moretodo}.
			
			\begin{figure}
				\centering
				\begin{tikzpicture}
					\begin{pgfonlayer}{nodelayer}
						\node [style=empty dot] (0) at (-9, -1) {$a$};
						\node [style=empty dot] (1) at (-8, -2) {};
						\node [style=empty dot] (2) at (-8, 0) {};
						\node [style=none] (3) at (-7, 0) {};
						\node [style=empty dot] (4) at (-4, 0) {$m$};
						\node [style=empty dot] (5) at (-3, -1) {$y$};
						\node [style=empty dot] (9) at (-5, 0) {$l$};
						\node [style=none] (12) at (-6, 0) {};
						\node [style=none] (13) at (-6.5, 0) {$\dots$};
						\node [style=empty dot] (14) at (-2, 3) {$b'$};
						\node [style=empty dot] (15) at (-2, 2) {$c'$};
						\node [style=empty dot] (16) at (-1, 1) {$m'$};
						\node [style=empty dot] (17) at (-3, 3) {$a'$};
						\node [style=empty dot] (18) at (0, 1) {};
						\node [style=empty dot] (19) at (-1, 0) {};
						\node [style=empty dot] (20) at (0, 0) {};
						\node [style=empty dot] (21) at (-2, -1) {};
						\node [style=empty dot] (22) at (-2, -2) {};
						\node [style=empty dot] (23) at (-9, -7) {$a$};
						\node [style=empty dot] (24) at (-8, -8) {};
						\node [style=empty dot] (25) at (-8, -6) {};
						\node [style=none] (26) at (-7, -6) {};
						\node [style=empty dot] (27) at (-4, -6) {$m$};
						\node [style=empty dot] (30) at (-5, -6) {$l$};
						\node [style=none] (31) at (-6, -6) {};
						\node [style=none] (32) at (-6.5, -6) {$\dots$};
						\node [style=empty dot] (33) at (-2, -4) {$b'$};
						\node [style=empty dot] (34) at (-2, -5) {$c'$};
						\node [style=empty dot] (36) at (-3, -4) {$a'$};
						\node [style=empty dot] (37) at (-0.5, -5.5) {};
						\node [style=empty dot] (38) at (-1.5, -6.5) {};
						\node [style=empty dot] (39) at (-0.5, -6.5) {};
						\node [style=empty dot] (40) at (-2, -7) {};
						\node [style=empty dot] (41) at (-2, -8) {};
						\node [style=none] (42) at (-11, -1) {$G_{i+1}$};
						\node [style=none] (43) at (-11, -7) {$G_{i+2}$};
						\node [style=empty dot] (44) at (-3, -7) {$y$};
					\end{pgfonlayer}
					\begin{pgfonlayer}{edgelayer}
						\draw (1) to (0);
						\draw (0) to (2);
						\draw (2) to (3.center);
						\draw (9) to (4);
						\draw (4) to (5);
						\draw (9) to (12.center);
						\draw (17) to (14);
						\draw (14) to (15);
						\draw (15) to (16);
						\draw (16) to (19);
						\draw (16) to (18);
						\draw (18) to (20);
						\draw (20) to (19);
						\draw (5) to (21);
						\draw (5) to (22);
						\draw (22) to (21);
						\draw (21) to (19);
						\draw (24) to (23);
						\draw (23) to (25);
						\draw (25) to (26.center);
						\draw (30) to (27);
						\draw (30) to (31.center);
						\draw (36) to (33);
						\draw (33) to (34);
						\draw (37) to (39);
						\draw (39) to (38);
						\draw (41) to (40);
						\draw (40) to (38);
						\draw (44) to (41);
						\draw (44) to (40);
						\draw (44) to (27);
					\end{pgfonlayer}
				\end{tikzpicture}
				\caption{\label{fig:moretodo} The second step of an application of rule \ref{it:gt}, in case B(ii).}
			\end{figure}
			
			This shows, essentially, that for the next steps in our construction of $W$ we either address a germ ending with $y$, removing $y$, or we address germs which are contained in $T$, and each time we do the latter action, we break up a circle among those which make up $T$. This turns the circle under consideration into a line segment, which is a connected component. Now, by Lemma \ref{lem:trees}, the number $k$ of circles which make up $T$ is bounded by the number of nodes of degree $2$ in $G_i$, which is $\leq 3i$. Therefore, by stage $4i$, we must address the germ $aEbE\dots y$ and remove $y$. Hence, no germ starts with $aEb$ in $G_{4i}$ and the connected component of $b$ in $G_{4i} \setminus \{a\}$ is a line segment.
		\end{proof}
		
		\begin{claimn}\label{claim:rightdegree}
			$\Delta(G \setminus W) =2$.
		\end{claimn}
		\begin{proof}
			Notice that $x$ has degree $3$ in $G \setminus W$ if and only if $\forall i \,\, x \in G_i \land \delta_{G_i}(x)=3$. This is a $\bPi{0}{1}$ condition, so by $\bSigma{0}{1}$-induction it follows that, if there are $x \in G \setminus W$ with degree $3$, then there is a least one, say $s$.
			
			This implies that there is $i_0$ such that for every $i>i_0$, $s$ is the least node of degree $3$ in $G_i$.
			
			We argue that it cannot be the case that, for cofinitely many $j > i_0$, we pick $w_j$ following rule \ref{it:pv} (so $w_j$ is a $P$-vertex in $G_j$). By contradiction, assume this were the case,  and let $h$ be such that for all $k \geq h$, $w_k$ is picked following rule \ref{it:pv}. Note that the deletion of a $P$-vertex from a graph $G_j$ reduces the number of nodes of degree $2$ (the two neighbors of $w_j$ of degree $2$ have degree $1$ in $G_{j+1}$, and compensate for the possibility that the third neighbor gets degree $2$ in $G_{j+1}$). Since there are only finitely many nodes in $G_h$ with degree $2$, this is a contradiction.
			
			So it must be the case that, infinitely often, $w_i$ is picked following some other rule. Clearly, for all $i > i_0$, $w_i$ cannot be the result of applying rule \ref{it:qv}, since rule \ref{it:qv} either deletes $s$ or deletes one of its neighbors, decreasing the degree of $s$.
			
			The only possibility is that, for every $i > i_0$, $w_i$ is picked using either Rule \ref{it:pv} or Rule \ref{it:germs}, and for infinitely many $j > i_0$, $w_j$ is picked using Rule \ref{it:germs}.
			
			Note that for every $j$, one can apply Rule \ref{it:pv} at most $3j$ consecutive times starting from $G_j$, and, by Claim \ref{claim:germs}, if $aEbE\dots$ is a germ in $G_j$ which is addressed at step $j$, there are no germs starting with $aEbE\dots$ in $G_{4j}$. This yields an explicit upper bound $k$ such that, in $G_k$, there are no germs starting with a node $a < s$ (it is not hard to see that $k=4^{4s}i_0$ works). This is a contradiction.
		\end{proof}
		
		\begin{claimn}\label{claimcompo}
			The connected components of $G \setminus W$ consist of line segments.
		\end{claimn}
		\begin{proof}
			By Claim \ref{claim:rightdegree}, $\Delta(G \setminus W)=2$, so its connected components can only be lines, line segments, or circles. By Claim \ref{claim:notree}, no connected component of $G \setminus W$ is a circle. So we only need to exclude lines. By contradiction, suppose that $P \subseteq G \setminus W$ is a connected component which is a line. Since each $G_i$ only has finitely many nodes of degree $2$, there must be some $a, b \in P$ and some $i_0$ such that, for every $i \geq i_0$, there is a germ in $G_i$ starting with $aEb$. This contradicts Claim \ref{claim:germs}.
		\end{proof}
		
		Now we can give an algorithm to find a $3$ coloring of $G$: given $v \in G$, we run the construction of $W$ until we either enumerate $v \in W$ (in which case we color $v$ with color $0$) or we enumerate a large enough finite piece $W_s \subseteq W$ such that the connected component of $v$ in $G \setminus W_s$ is a line segment. Claim \ref{claim:independent} guarantees that the component of $v$ in $G \setminus W_s$ coincides with the component of $v$ in $G \setminus W$. Hence, we can just $2$-color (with colors $1$ and $2$) this line segment starting from its least endpoint.
	\end{proof}
	
	\subsection{$\rca$ proves $\forall d \geq 3 \,\, \btb[d]$.}
	
	We are now ready to prove, by internal induction on $d$, that $\rca$ proves the restriction of Brooks' Theorem to bounded graphs with $d \geq 3$. The key ingredients of this induction are already present in the proof of Theorem \ref{thm:bt3}.
	
	\begin{proposition}\label{prop:ddown}
		There are a $\lSigma{0}{1}$ formula $\varphi$ and a $\lPi{0}{1}$ formula $\psi$ such that, for every $d \geq 4$ and every $d$-regular $G=(V,E)$ which does not induce any $d+1$ clique, we have
		\[\rca \vdash \forall n \,\, (\varphi(d,G,n) \leftrightarrow \psi(d, G, n)),\]
		and the set $H$ which is defined by $\varphi$ induces a subgraph $H \subseteq G$ with $\Delta(H)=d-1$ which does not induce any $d$-clique.
	\end{proposition}
	
	In particular, by formalized computability theory within $\rca$, we talk about a uniform reduction $\Phi$ such that, if $G$ is $d$-regular and satisfies the hypothesis for Brooks' Theorem, then $\Phi(d,G)$ is the characteristic function of a bounded subgraph $H \subseteq G$, satisfying the hypothesis of $\btb{d-1}$.
	
	\begin{proof}[Proof of Proposition \ref{prop:ddown}]
		We provide an algorithmic procedure to enumerate a set $W=(w_i : i \in \N)$ of vertices of $G$. This procedure is similar to that of the proof of Theorem \ref{thm:bt3}. This set $W$ is prima facie only computably enumerable in $G$ (i.e.\ defined by a $\lSigma{0}{1}(G)$ formula), but we will again have:
		\begin{itemize}
			\item $W$ is an independent subset of $G$, and
			\item $\Delta(G \setminus W)=d-1$,
		\end{itemize}
		so reasoning as in the proof of Theorem \ref{thm:bt3} we see that $W$ is actually computable ($\lDelta{0}{1}(G)$-definable) uniformly in $G$ and $d$.
		The procedure to enumerate $(w_i : i \in \N)$ is based on the notions of $P$-vertices and $Q$-vertices, adapted to the degree $d$. In this context, a $P$-vertex of a graph $G$ is a vertex $v \in G$ which belongs to an induced $d$-clique $K$, has degree $d$, and is the only vertex of $K$ with degree $d$. A $Q$-vertex is a vertex $v \in G$ which belongs to an induced $d$-clique $K$ with at least $2$ and at most $d-1$ vertices of degree $d$. Lastly, $d$-trees are objects which are built out of $d$-cliques in the same way that circle-trees are built out of circles.
		
		The enumeration of $(w_i : i \in \N)$ proceeds as follows: given $G_i=G \setminus \{w_j : j < i\}$, we pick $w_i$ to be the least $P$-vertex of $G_i$, if there are any, otherwise we look at the least vertex $v \in G_i$ with $\delta_{G_i}(v)=d$, and we pick the least $Q$-vertex among the neighbors of $v$, if there are any, otherwise, we let $w_i=v$.
		
		As the reader will be able to see, this is a simplified version of the algorithm provided in the proof of Theorem \ref{thm:bt3} (we omit the analogue of case \ref{it:germs} in \ref{thm:bt3}).
		
		Minor modifications of the argument presented in the aforementioned proof (e.g.\ those necessary to deal with trees built out of $d$-cliques rather than circles) are sufficient to show that, letting $H= G \setminus W$
		\begin{itemize}
			\item the procedure just described is effective,
			\item there are no $d$ cliques in $H$,
			\item $W$ is an independent set,
			\item $\Delta(H)=d-1$.
		\end{itemize}
		Indeed, the these statements correspond to Claims \ref{claim:effectiveness}, \ref{claim:notree}, \ref{claim:independent}, and \ref{claim:rightdegree} of that proof.
	\end{proof}
	\begin{remark}\label{rem:comb}
		Combining the algorithm giving the function $G \mapsto H$ with the algorithm of Lemma \ref{lem:boundedregular} which takes a bounded graph $H$ with $\Delta(H)=d-1$ and embeds it into a $(d-1)$-regular graph $H'$, we obtain a $(d-1)$-regular graph $H'$ which satisfies the hypotheses of $\btb{d-1}$ and includes $H$ as an induced subgraph.
	\end{remark}
	
	\begin{definition}
		Assuming that $G=(V,E,f)$ is a bounded graph with $\Delta(G)=d \geq 2$, we denote by $r(d,G)$ the encoding of the result of the proof of Lemma \ref{lem:boundedregular}: a $d$-regular graph $G'=(V', E')$, a set $U \subseteq V'$ and a bijection $j(d,G) \colon V \rightarrow U$ which is a graph isomorphism from $G$ to the induced subgraph of $G'$ determined by $U$. By $r'(d,G)$ we denote the graph $G'$ described above.
	\end{definition}
	
	Given a graph $G$ with $\Delta(G)=d \geq 3$ which does not induce any $d+1$-cliques, we exploit Proposition \ref{prop:ddown} and Remark \ref{rem:comb} to define by recursion a sequence of sets $(K(i) : i < d-3)$ such that each $K(i)$ is an independent subset of $G$, $G_3=G \setminus \bigcup_{i < d-3} K(i)$ is a bounded graph with $\Delta(G_3)=3$, and $G_3$ does not induce any $4$-clique.
	
	\begin{definition}\label{def:color}
		Let $d \geq 3$ and let $G$ be a bounded graph with $\Delta(G)=d$ which does not induce any $d+1$-cliques. We define a sequence of graphs as follows: we let $H_{0}=G$, and for every $i < d-3$, we let $H_{i+1}=\Phi(d-i, r'(d-i,H_i))$.
	\end{definition}
	
	We show that Definition \ref{def:color} is a legitimate $\rca$ definition by induction.
	
	\begin{lemma}\label{lem:partialcoloring}
		Let $d \geq 3$ and let $G$ be a graph which satisfies the hypotheses of $\btb{d}$. For every $i \leq d-3$ the graph $H_{i}$ satisfies the hypothesis of $\btb{d-i}$.
	\end{lemma}
	\begin{proof}
		Our claims can be formalized as $\bPi{0}{1}$ formulas, so they can be proved by induction. The inductive step follows from Lemma \ref{lem:boundedregular} and Proposition \ref{prop:ddown}.
	\end{proof}
	
	Given a graph $G$ which satisfies the hypothesis of $\btb{d}$, we can ``follow'' the nodes of $G$ in the sequence of graphs $(H_i : i \leq d-3)$ using the injections $j(d-i, H_{i})$ of Definition \ref{def:color}, and thus we can define a sequence of sets $(K(i) : i < d-3)$ with $K(i) \subseteq V$ as $v \in K(i)$ if and only if $i$ is the least such that the image of $v$ via the composition of these injections does not belong to $H_{i+1}$. Again by induction, using Proposition \ref{prop:ddown}, we know that $(K(i) : i < d-3)$ is a sequence of independent subgraphs of $G$, and $G_3= G \setminus \bigcup_{i <d-3}K(i)$ is a bounded graph with $\Delta(G_3)=3$ which does not induce any $4$-cliques.
	
	\begin{theorem}\label{thm:btd}
		$\rca$ proves $\forall d \geq 3 \,\, \btb d$.
	\end{theorem}
	\begin{proof}
		Let $G$ be a graph satisfying the hypotheses of $\btb{d}$. If $\Delta(G)=3$, then $G$ is $3$-colorable by Theorem \ref{thm:bt3}. Otherwise, consider the sequence $(K(i) : i < d-3)$ and the bounded graph $G_3$ defined above. By Theorem \ref{thm:bt3}, there is a proper $3$-coloring of $G_3$. Gluing together this coloring with the coloring naturally given by the sets $(K(i) : i < d-3)$ yields a proper $d$-coloring of $G$.
	\end{proof}

	\appendix
	
	\section{On Schmerl's effective proof of Brooks' Theorem}
	
	A search through the literature reveals two different proofs of an effective version of Brooks' Theorem, one by Tverberg \cite{Tverberg84}, and one by Schmerl \cite{Schmerl82}. Indeed, Tverberg's article is presented as a response to Schmerl's, giving a simplified proof of the result. Careful scrutiny revealed an error in Schmerl's proof: there is a counterexample to the statement of his Lemma 7.
	
	Schmerl's proof is structured in a way similar to Tverberg's: there is an induction on the number $\Delta(G)$, and the base case where $\Delta(G)=3$ is proved separately. Lemma 7 is used crucially in proving Brooks' theorem for $\Delta(G)=3$. We do not know whether there is a fix to this problem. Here we present Schmerl's Lemma, and our counterexample.
	
	\begin{definition}
		Let $G=(V,E)$ be any graph. We define the graph $\dlow(G)$ as the induced subgraph of $G$ with vertex set $V'=\{v \in V : \delta_G(v) \leq 2\}$
	\end{definition}
	\begin{definition}
		A finite graph $G=(V,E)$ has property A if $\Delta(G)=3$ and any $3$-coloring $c \colon \dlow(G) \rightarrow 3$ can be extended to a $3$-coloring of $G$.
	\end{definition}
	\begin{lemma}\label{sch7}
		Let $H$ be a finite graph such that $\Delta(H) \leq 3$ and let $Z_0 \subseteq H$ be such that there is some $h \in \N$ such that $N_h(Z_0)=H$. For all $i < h$ let $Z_{i+1}=N_{i+1}(Z_0) \setminus N_i(Z_0)$. Assume that the following hold:
		\begin{enumerate}
			\item $h$ is sufficiently large,
			\item for every $i < h$, $|\{x \in Z_i : N_1(x) \cap Z_{i+1} \neq \emptyset\}| \leq 2$,
			\item $1 \leq |\dlow(H) \cap Z_h| \leq 2$,
			\item $Z_0 \subseteq \dlow(H) \subseteq Z_0\cup Z_h$.
		\end{enumerate}
		Then $H$ has property A.\footnote{The original paper presents $Z_0 \cap Z_h$ on item $(4)$, and the first mention of $Z_0$ is not capitalized. These are obvious typos.}
	\end{lemma}
	A graph satisfying the hypotheses of Lemma \ref{sch7} is naturally divided in $h+1$ ``levels'' (the $Z_i$'s). We refer to $Z_0$ as the first level of such a graph. Schmerl claims that if $h \geq 6$, the Lemma holds.
	
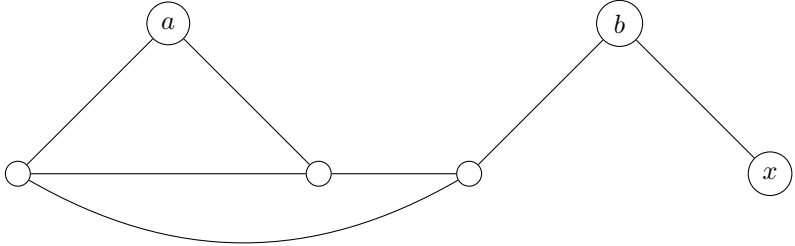
	\begin{figure}
		\centering
		\begin{tikzpicture}
			\begin{pgfonlayer}{nodelayer}
				\node [style=empty dot] (0) at (-5, 3) {$a$};
				\node [style=empty dot] (1) at (-7, 1) {};
				\node [style=empty dot] (2) at (-3, 1) {};
				\node [style=empty dot] (3) at (1, 3) {$b$};
				\node [style=empty dot] (4) at (-1, 1) {};
				\node [style=empty dot] (5) at (3, 1) {$x$};
			\end{pgfonlayer}
			\begin{pgfonlayer}{edgelayer}
				\draw (1) to (0);
				\draw (0) to (2);
				\draw (1) to (2);
				\draw (2) to (4);
				\draw [bend left] (4) to (1);
				\draw (4) to (3);
				\draw (3) to (5);
			\end{pgfonlayer}
		\end{tikzpicture}
		\caption{\label{fig:counterexample} Graph $C$, which we use to build counterexamples to Lemma \ref{sch7}.}
	\end{figure}

	Inspection of the conditions imposed on the graph $H$ reveals that one can have $Z_0$ consist of a single node, and $h \in \N$ arbitrarily large. Now consider the graph $C$ in Figure \ref{fig:counterexample}.
	Given any graph $G$ which satisfies the hypotheses of Lemma \ref{sch7} and such that its ``first level'' $Z_0$ consists of a single node $g$, we can obtain a graph $CG$ by considering the disjoint union $C \cup G$ quotiented by the equivalence relation identifying $g$ and $x$. Then the graph $CG$ also satisfies the hypotheses of Lemma \ref{sch7}; the first level of $CG$ consists of nodes $\{a,b\}$, and, if $G$ consists of $h$ levels, $CG$ consists of $h+1$ levels. It is immediate to see that we can $3$-color $\dlow(CG)$ assigning the same color to $a$ and $b$, and this coloring cannot be extended to the entirety of $CG$ (it cannot even be extended to $C$). This shows that there is no $h$ such that Lemma \ref{sch7} holds for graphs consisting of more than $h$ levels.

\bibliographystyle{alpha}

\end{document}